\documentclass[10pt,oneside]{amsbook}
\usepackage{amsmath}
\usepackage{breqn}
\usepackage{multicol}
\usepackage{tikz}
\usepackage{pgfplots}
\pgfplotsset{compat=1.8}
\usepgfplotslibrary{fillbetween}
\usepackage{collectbox}

\usepackage{array}

\usepackage[hidelinks]{hyperref}

\usepackage{amssymb}

\usepackage{skak}
\usepackage{graphicx}
\usepackage{tikz}

\usepackage{mathtools}

\DeclarePairedDelimiter\floor{\lfloor}{\rfloor}

\usepackage{stackengine}
\newlength\Lobj
\newsavebox{\newobj}

\newcommand\LoF[2]{%
  \def\obj{$#1$~$#2$\,}%
  \setlength\Lobj{\widthof{\obj}}%
  \sbox\newobj{\stackon[3pt]{\obj}{\rule{\Lobj}{.3pt}}}%
  \usebox{\newobj}%
  \rule{.3pt}{\ht\newobj}%
  }

\numberwithin{equation}{chapter}
\numberwithin{figure}{chapter}
\theoremstyle{plain} % This is the default style.
\newtheorem{thm}{Theorem}
\begin{document}
\frontmatter % Pages here will be numbered with Roman numerals.

\vspace*{1 cm} 
% LaTeX will remove blank space at the start or end of a page,
%   but the asterisk prevents this.

\begin{center}
 
 \LoF{\vphantom{@}}{\hphantom{.}} \LoF{\vphantom{@}}{\hphantom{.}} = \LoF{\vphantom{@}}{\hphantom{.}}\\
 
 \medskip
 
 \LoF{\LoF{\vphantom{@}}{\hphantom{.}}}{\!}  =  \hphantom{\LoF{\vphantom{@}}{\hphantom{.}} \LoF{\vphantom{@}}{\hphantom{.}}} \\

\end{center}

\bigskip\bigskip

\centerline{\Huge Laws Of Form and the Riemann Hypothesis}
\medskip
\bigskip\bigskip
\bigskip\bigskip

\centerline{\LARGE J. M. Flagg}
\medskip
\medskip
%\centerline{\Large (George Spencer-Brown)}
\medskip
\medskip
\centerline{\LARGE Louis H. Kauffman}
\begin{center}
Department of Mathematics, Statistics, and Computer Science,\\
851 South Morgan Street,\\
University of Illinois at Chicago,\\
Chicago, IL 60607-7045\\
and\\
Department of Mechanics and Mathematics,\\
Novosibirsk State University\\
Novosibirsk,\\
Russia\\
kauffman@uic.edu
\end{center}
\medskip
\medskip
\centerline{\LARGE Divyamaan Sahoo}
\begin{center}
Sound Department,\\
School of the Art Institute of Chicago,\\
112 S Michigan Ave, Suite 512,\\
Chicago, IL 60603\\
dsahoo@artic.edu
\end{center}

%\medskip
%\centerline{\Large Chicago, Illinois}
%\medskip
%\centerline{\Large \textit{June, 2020}}

%\begin{abstract}
%Abstract
%\end{abstract}

\newpage

\newpage

\thispagestyle{empty}

\section*{0. Abstract}

\noindent This paper is an exposition of the Riemann Hypothesis from a historical point of view, and an investigation of part of the work George Spencer-Brown has done on the Riemann Hypothesis. We concentrate on Appendix 9 of the revised sixth edition of Laws of Form (2016) and plan to discuss Appendix 7 and 8 in papers to follow.   Spencer-Brown's approach to number theory is placed in context in this paper. His approach to the Riemann Hypothesis is described here in the context of the evolution of the problem, and in particular in relation to the reformulation of RH in terms of the behaviour of the cumulative M\"{o}bius Function $M(n)$. In this light, we know that to prove RH is to show that $M(n)$ is of the order $O(\sqrt{n} n^{\epsilon})$. On the other hand, Spencer-Brown has given a profound formula for finding the M\"{o}bius function in terms of its previous values. This formula has the form
$$ \mu(n) =  - \Bigg(\sum_{d=1}^{n-1} \mu(d)\floor*{\frac{n}{d}}\Bigg),$$
where for $d=1$, we make the Spencer-Brown adjustment, rectifying the Floor function by replacing $\big[\frac{n}{d}\big]$ by $\big[\frac{n-1}{d}\big]$.

\noindent The key point about this formula is that it is a negative feedback for the next value of the M\"{o}bius function with respect to all of its previous values. The formula shows that the M\"{o}bius function is far from being random and indeed that any run of strictly positive values or any run of strictly negative values will be turned around by the negative feedback. This means that a cumulative sum of the values of the M\"{o}bius function will exhibit highly controlled oscillation about zero, with widening swings as the value of $n$ in $M(n)$ grows. The recursion implicit in the Spencer-Brown/Legendre formula makes $M(n)$ self-regulating and it is this self-regulation that drives its asymptotic behaviour. Spencer-Brown's fundamental observation is that the M\"{o}bius function $\mu (n)$, while not random, behaves in a more constrained way than a random series due to the feedback inherent in its dependence on its previous values. To quote Spencer-Brown,  ``$\cdots$ the successive nonzero values of $\mu(d)$ behave like the falls of a magic coin that, from the moment it is struck, remembers exactly how many times it has fallen with one side or the other uppermost, and whenever one side exceeds the other, biases itself towards the other until the excess is eliminated" (page 211, Appendix 9 \cite{SBG}). This means that the asymptotics for a random series still apply to the M\"{o}bius function and so the Riemann Hypothesis (RH) follows.

\bigskip

\begin{center}
   \textit{In dedication to John Albert Ewell III (February 28, 1928 - July 21, 2007)}
\end{center}

\newpage

\mainmatter % Pages here will be numbered 1, 2, 3, and so on .

\thispagestyle{empty}

\section{Introduction to the Riemann Hypothesis}

\noindent The Riemann Hypothesis (RH) is a conjecture first made by Bernhard Riemann in \textit{\"{Uber} die Anzahl der Primzahlen unter einer gegebenen Gr\"{o}sse}, ``On the Number of Primes Less Than a Given Magnitude", 1859 \cite{RBU}. The conjecture is about the ``zeros" of the zeta function, $\zeta$, whose domain is  the complex numbers, $s \in \mathbb{C}$:\

$$\zeta (s) = \sum_{1}^{\infty}\frac{1}{n^s} =  \frac{1}{1^s} + \frac{1}{2^s} + \frac{1}{3^s} + ... ,$$\

\noindent where $a, b \in \mathbb{R}$, and $i = \sqrt{-1}$.

\noindent By convention, the trivial zeros of $\zeta$ are the integers $-2, -4, -6, -8, ...$ . These are understood to be zeros via analytic continuation\footnote{In Section 1.4, ``The Function $\zeta(s)$", \emph{Riemann's Zeta Function} \cite{EHM}, H. M. Edwards writes: ``The view of analytic continuation in terms of chains of disks and power series convergent in each disk descends from Weierstrass and is quite antithetical to Riemann's basic philosophy that analytic functions should be dealt with \textit{globally}, not locally in terms of power series". Edwards notes that ``Riemann does not speak of the ``analytic continuation" of the function $\Sigma n^{-s}$ beyond the halfplane Re $s>1$, but speaks rather of finding a formula for which it ``remains valid for all of $s$". In \textit{Ueber die Anzahl der Primzahlen unter einer gegebenen Gr\"{o}sse} \cite{RBU}, Riemann indeed writes: ``The function of the complex variable $s$ which is represented by these two expressions (\textit{the two halves of the Euler product identity}), wherever they converge, I denote by $\zeta (s)$. Both expressions converge only when the real part of $s$ is greater than $1$; at the same time an expression for the function can easily be found which always remains valid" (transl. Wilkins, 1998 \cite{RBU}).}. The remaining zeros of $\zeta$ are called the non-trivial zeros. Riemann conjectured that the non-trivial zeros of zeta are always of the form $s = \frac{1}{2} + bi$, i.e., that the ``real part" $a$ of the non-trivial zero $s$, is  $\frac{1}{2}$.

\noindent Euler first studied the zeta function in the domain of the real numbers\footnote{In Section 1.2 ``The Euler Product", \textit{Riemann's zeta function} \cite{EHM}, H. M. Edwards writes, ``it seems certain that Riemann's use of the Euler product formula was influenced by Dirichlet" and ``Dirichlet, unlike Euler, used the formula (1) with $s$ as a real variable and, also unlike Euler, he proved rigorously that (1) is true for all real $s>1$" (Page 7, Edwards). Edwards provides Dirichlet's remarks in a footnote: ``Since the terms $p^{-s}$ are all positive, there is nothing subtle or difficult about this proof -- it is essentially a reordering of absolutely convergent series -- but it has the important effect of transforming (1) from a formal identity true for various values of $s$ to an analytical formula true for all real $s>1$" (Dirichlet, page 7, footnote by Edwards). We know, of course, that Riemann went even further, broadening the scope of $\zeta$ to include all complex $s$, $s \neq 1$. Here, formula $(1)$ is the Euler product formula provided above.}, $\mathbb{R}$. In \textit{Variae observationes circa series infinitas} (1737) \cite{E}, he derived a fundamental relationship between the zeta function and the prime numbers: %Consider,\\

$$\zeta (s) = \sum_{n}^{} \frac{1}{n^s} = \frac{1}{1^s}+\frac{1}{2^s}+\frac{1}{3^s}+... =  \frac{1}{(1-\frac{1}{2^s})(1-\frac{1}{3^s})(1-\frac{1}{5^s})...\hphantom{.}} = \prod_{p}^{} \frac{1}{1-p^{-s}}.$$

\noindent Riemann extends $\zeta$ to the complex domain by using the factorial function, notated $\Pi(s)$ due to Gauss in \textit{Disquisitiones generales circa seriem infinitam} (1813) \cite{G}. Originally attributed to Euler in \textit{De progressionibus transcendentibus seu quarum termini generales algebraice dari nequeunt}(1730) \cite{ED}\footnote{Euler wrote it as follows: $y = e^{-x}$ and $n! = \int_{0}^{1} (log \frac{1}{y})^n dy$}, the factorial function for $n \in \mathbb{N}$,
$$ n! = \int_{0}^{\infty} e^{-x}x^n dx$$
can be extended to the real numbers $> -1.$ Due to Gauss \cite{G} we have
$$ \Pi (s) = \int_{0}^{\infty} e^{-x}x^s dx, \quad \quad \quad (s > -1).$$

\noindent Abel's \textit{Solution de quelques problemes a l’aide d‘integrales d\'{e}finies} (1823) \cite{ANH} and Chebyshev's \textit{Sur la fonction qui d\'{e}termine la totalit\'{e} des nombres premiers
inferieurs \`{a} une limite donn\'{e}e} (1852) \cite{CH} give us the following identity for real values of $s$, $s>1$:

$$ \zeta (s) = \sum_{1}^{\infty} \frac{1}{n^s} = \frac{\Pi (-s)}{2\pi i} \int_{+\infty}^{+\infty} \frac{(-x)^s}{e^x - 1}.\frac{dx}{x}. $$
\newpage
\noindent The contour integral in $\frac{\Pi (-s)}{2\pi i} \int_{+\infty}^{+\infty} \frac{(-x)^s}{e^x - 1}.\frac{dx}{x}$ clearly converges for all values of $s$, real or complex, since $e^x$ grows faster than $x^s$ as $x \to \infty$. Since $\Pi (s)$ is an analytic function of the complex variable $s$, which has simple poles at $s=-1, -2, -3, ...$, $\Pi (-s)$ has poles at $s = 1, 2, 3, ...$\hphantom{.}. We know $\zeta(s)$ converges for $s = 2, 3, 4, ...$ (\`{a} la Euler et al), so by Cauchy's theorem, the integral must have a zero which cancels the pole of $\Pi (-s)$ at $s = 2, 3, 4, ...$\hphantom{.}. Thus, $\zeta (s) = \sum_{1}^{\infty} \frac{1}{n^s} = \frac{\Pi (-s)}{2\pi i} \int_{+\infty}^{+\infty} \frac{(-x)^s}{e^x - 1}.\frac{dx}{x}$ is analytic at all points of the complex $s$-plane except for a simple pole at $s = 1$ (Edwards, Section 1.4, ``The Function $\zeta(s)$" \cite{EHM}).

\noindent On page 2-3 of \textit{Ueber die Anzahl der Primzahlen unter einer gegebenen Gr\"{o}sse} (transl. Wilkins, 1998) \cite{RBU}, Riemann continues, ``thus a relation between $\zeta(s)$ and $\zeta(1-s)$, which through the use of known properties of the function $\Pi$, may be expressed as follows" and that ``remains unchanged when $s$ is replaced by $1-s$":

$$ \Pi \bigg( \frac{s}{2} - 1 \bigg) \pi^{-\frac{s}{2}} \zeta (s) = \Pi \bigg( \frac{1-s}{2} - 1 \bigg) \pi^{-\frac{1-s}{2}} \zeta (1-s).$$
\\
Riemann considers this symmetrical statement as the natural statement of the functional equation. The function $\Pi \big(\frac{s}{2}-1\big)\pi^{-\frac{s}{2}}\zeta (s)$ in symmetrical form, has poles at $s=0$ and $s=1$. Riemann then introduces $\xi$, an analytic function defined for all values of $s$, by multiplying the simpler left hand side of the functional equation of $\zeta$ by $\frac{s(s-1)}{2}$ to get\footnote{Adopting Landau's notational convention for $\xi$, Edwards remarks ``Actually Riemann uses the letter $\xi$ to denote the function which it is now customary to
denote by $\Xi$, namely, the function $\Xi(t) = \xi(\frac{1}{2} + it)$, where $\xi$ is defined as above." (Edwards, footnote on page 16 \cite{EHM})}

$$ \xi(s) = \Pi\bigg(\frac{s}{2}\bigg)(s-1)\pi^{-\frac{s}{2}} \zeta (s)$$
\\
Given a convergent infinite product can only be zero if one of its factors is zero, $\zeta(s) = \Pi_p (\frac{1}{1-p^{-s}})^{-1}$ cannot have zeros for Re $s>1$. Since the factors other than $\zeta(s)$ in $\xi(s) = \prod(\frac{s}{2})(s-1)\pi^{-\frac{s}{2}}\zeta(s)$ only have the simple zero at $s=1$, we know that none of the roots $\rho$ of $\xi(\rho) = 0$ lie in Re $s>1$. Due to the functional equation, we know $\rho$ is a root \emph{iff} $1-\rho$ is a root, so we also know that none of the roots $\rho$ of $\xi(\rho) = 0$ lie in the half-plane Re $s<0$. Hence, all the roots $\rho$ of $\xi(\rho) = 0$ lie in the strip $0 \leq$ Re $\rho \leq 1$ (Edwards, Section 1.9 ``The Roots $\rho$ of $\xi$" \cite{EHM}).

\newpage

\noindent The following product representation of $\xi(s)$ was proved by Hadamard in \textit{\'{E}tude sur les Propri\'{e}t\'{e}s des Fonctions Enti\`{e}res et en Particulier d'une Fonction Consid\'{e}r\'{e}e par Riemann} (1893) \cite{HJ}\footnote{Since the function $\log \xi(s)$ has logarithmic singularities at the roots $\rho$ of $\xi(s)$ and no other singularities, it has the same singularities as the formal sum $\sum_{\rho} \log(1-\frac{s}{\rho})$. If this sum converges and is well behaved with respect to $\log \xi(s)$ near $\infty$, then $\sum_{\rho} \log(1-\frac{s}{\rho})$ differs from $\log \xi(s)$ by at most an additive constant. Setting $s = 0$ gives $\log \xi(0)$, which upon exponentiating gives us the identity (H. M. Edwards, Section 1.10 ``The Product Representation of $\xi(s)$" \cite{EHM}).\\
\indent For any fixed $s$, the ambiguity in the imaginary part of $\log\big[1-\frac{s}{\rho}\big]$ disappears for large $\rho$. Hence, the sum $\sum_{p}^{}\log\big(1-\frac{s}{\rho}\big)$ is defined except for a finite multiple of $2\pi i$ which drops out when one exponentiates the $\xi$ identity above. Furthermore, ``one can ignore the imaginary parts altogether", the real parts of the terms of $\sum_{p}^{}\log\big(1-\frac{s}{\rho}\big)$ ``are unambiguously defined and their sum is a harmonic function which differs from $\text{Re} \log \xi(s)$ by a harmonic function without singularities, and if this difference function can be shown to be constant, it will follow that its harmonic conjugate is constant also" (Edwards, Section 1.10 ``The Product Representation of $\xi(s)$" \cite{EHM}).\\
\indent There are two problems associated with $\sum_{\rho} \log(1-\frac{s}{\rho})$, the first being the ``determination of the imaginary parts of the logarithms it contains", and the second associated with its convergence, for the sum is ``in fact a conditionally convergent sum, and the \textit{order} of the series must be specified in order for the sum to be well determined", i.e., it suffices merely to stipulate that each term be paired with its ``twin" $\rho \leftrightarrow 1-\rho$. (Edwards, Section 1.10 ``The Product Representation of $\xi(s)$" \cite{EHM})}:
$$ \xi(s) = \xi(0) \prod_{\rho}^{}\bigg( 1 - \frac{s}{\rho}\bigg),$$
\\
\noindent where $\rho$ ``ranges over" the roots of the equation $\xi(\rho)=0$ and multiple roots are to be counted with multiplicities. Recall that a polynomial $p(s)$ can be expanded as a finite product: 
$$p(s) = p(0)\prod_{\rho}^{}\bigg[1-\frac{s}{\rho}\bigg],$$

\noindent where $\rho$ ranges over the roots of the equation $p(\rho)=0$. Thus, $\xi(s)$ is ``\textit{like a polynomial of infinite degree\footnote{The identity $\xi(s) = \sum_{n=0}^{\infty}a_{2n}\big(s-\frac{1}{2}\big)^{2n}$ also suggests that $\xi(s)$ is a polynomial of infinite degree, and a finite number of terms of this series representation of $\xi(s)$, presented in Riemann's second proof of the functional equation, gives a very good approximation in any finite part of the plane.}}" (H. M. Edwards, Section 1.10, ``The Product Representation of $\xi(s)$" \cite{EHM}).

\noindent The $J$ function\footnote{``Riemann denoted this function $f(x)$, and most other writers denote it $\Pi(x)$. Since $f(x)$ now is commonly used to denote a generic function, I have taken the liberty of introducing a new notation $J(x)$ for this function" (Edwards, Section 1.11, footnote, page 22 \cite{EHM}).\\
``Let $F(x)$ be equal to this number when $x$ is not exactly equal to a prime number; but let it be greater by $\frac{1}{2}$ when $x$ is a prime number, so that, for any $x$ at which there is a jump in the value in $F(x), \hphantom{...} F(x) = \frac{F(x+0)+F(x-0)}{2}$. If in the identity $\log \zeta (s) = - \sum \log(1-p^{-s})=\sum p^{-s} + \frac{1}{2}p^{-2s}+\frac{1}{3}p^{-3s}+...$ one now replaces 
$$p^{-s} \hphantom{...} \text{by } s\int_{p}^{\infty}x^{-s-1}ds, \hphantom{...} p^{-2s}\hphantom{...} \text{by } s\int_{p^2}^{\infty}x^{-s-1}ds, \hphantom{...} ... \hphantom{...},\hphantom{...} \text{one obtains}$$
$$\frac{\log\zeta(s)}{s} = \int_{1}^{\infty}f(x)x^{-s-1}dx, \text{ if one denotes } F(x) + \frac{1}{2}F(x^{\frac{1}{2}})+\frac{1}{3}F(x^{\frac{1}{3}})+... \hphantom{...} \text{by } f(x).$$

This equation is valid for each complex value $a
+bi$ of $s$ for which $a>1$" (page 5-6, \textit{Ueber die Anzahl der Primzahlen unter einer gegebenen Gr\"{o}sse}, Riemann, transl. Wilkins, 1998 \cite{RBU}).} can be introduced via Stieltjes integrals. Using the identity\\
$\log (1-x) = -x -\frac{1}{2}x^2 - \frac{1}{3}x^3 ...$, the Euler product formula can be reformulated:
$$\log \zeta(s) = \sum_{p}^{}\bigg[\sum_{n}^{}\bigg(\frac{1}{n}\bigg)p^{-ns}\bigg] = \sum_{p}^{}\sum_{n}^{}\bigg(\frac{1}{n}\bigg)p^{-ns} = \int_{0}^{\infty} x^{-s} d J(x)$$
\noindent ``where $J(x)$ is the function which begins at $0$ for $x=0$ and increases by a jump of $1$ at primes $p$, by a jump of $\frac{1}{2}$ at prime squares $p^2$, by a jump of $\frac{1}{3}$ at prime cubes, etc" (Edwards, Section 1.11, ``The Connection Between $\zeta(s)$ And Primes" \cite{EHM}). Assuming the notation $\pi$ (attributed to Landau, \textit{Handbuch}, 1909 \cite{LE}) for the number of primes less than a given magnitude, we have the following two identities: 
$$ J(x) = \pi(x) + \frac{1}{2}\pi(x^{\frac{1}{2}})+\frac{1}{3}\pi(x^{\frac{1}{3}})+ ... +\frac{1}{n}\pi(x^{\frac{1}{n}}) + ...$$
$$ \pi(x) = \sum_{n=1}^{\infty}\frac{\mu(n)}{n} J(x^{\frac{1}{n}}) = J(x) - \frac{1}{2}J(x^{\frac{1}{2}}) - \frac{1}{3}J(x^{\frac{1}{3}}) - \frac{1}{5}J(x^{\frac{1}{5}}) + \frac{1}{6}J(x^{\frac{1}{6}})- ...$$
$$\text{It is of special interest to note that }\frac{1}{\zeta(s)} =\prod_{p}\bigg( 1 - \frac{1}{p^s}\bigg).$$
In the Euler product of $\big(1-\frac{1}{p^s}\big)$, we see the general term is $(-1)$ raised to the number of factors, a product of distinct prime factors. This leads to the definition $\mu(n) = (-1)^k$, where $k$ is the number of distinct prime factors of $n$, and $n$ has no higher order of factors, with $\mu(n) = 0$ otherwise\footnote{The classical M\"{o}bius function found in \emph{Uber ein besondere Art von Umkehrung der Reihen} (1832) \cite{MAF} is defined as follows: $g(m) = \sum_{n \rvert d} f(n)$ and $f(m) = \sum_{n \rvert d} g(n) \mu\big(\frac{m}{n}\big)$}.

$$\text{Thus, we have }\frac{1}{\zeta(s)} = \sum_{n=1}^{\infty}\frac{\mu(n)}{n^s}.$$

%\footnote{The classical M\"{o}bius function $\mu(n)$, found in \emph{Uber ein besondere Art von Umkehrung der Reihen} (1832) \cite{MAF} is: $$g(m) = \sum_{n\rvert d}^{} f(n) \text{and } f(m) = \sum_{n\rvert d}^{} g(n) \mu\big(\frac{m}{n}\big)}$$}

\newpage

\noindent In his second proof of the functional equation, Riemann performs the change of variable $x = n^2\pi x$ in Euler's integral for $\Pi(\frac{s}{2}-1)$. For Re $s>1$, he arrives at\\
$\Pi\big(\frac{s}{2}-1\big)\pi^{-\frac{s}{2}}\frac{1}{n^s} = \int_{0}^{\infty} e^{-n^2\pi x}x^{\frac{s}{2}}\frac{dx}{x}$. Summation over $n$ (for Re $s>1$) gives:
$$ \Pi\bigg(\frac{s}{2}-1\bigg)\pi^{-\frac{s}{2}}\zeta(s) = \int_{0}^{\infty} \psi (x) x^{\frac{s}{2}} \frac{dx}{x}$$
\noindent where $\psi(x) = \sum_{n=1}^{\infty}e^{-n^2\pi x}$. In order to prove that the function above is unchanged by the substitution $s = 1-s$, Riemann uses the form of the functional equation of Jacobi's theta function, referring to Section 65 of Jacobi's \textit{Fundamenta Nova Theoriae Functionum Ellipticarum}, in particular ``Suites des notices sur les fonctions elliptiques" (1828) \cite{J}. The following form is attributed by Jacobi to Poisson:
$$ \frac{1 + 2 \psi(x)}{1 + 2 \psi (\frac{1}{x})} = \frac{1}{\sqrt{x}}$$
\noindent We see that $\psi(x)$ approaches zero very rapidly as $x \to \infty$. This shows in particular that $\psi(x)$ is like $\frac{1}{2}(x^{-\frac{1}{2}}-1)$ for $x$ near zero. Hence $\int_{0}^{\infty}\psi(x)x^{\frac{s}{2}}\frac{dx}{x}$ is convergent for $s>1$. Thus\footnote{Riemann reformulates $\int_{0}^{\infty}\psi(x)x^{\frac{s}{2}}\frac{dx}{x} = \int_{1}^{\infty}\psi(x)\big[x^{\frac{s}{2}}+x^{\frac{1-s}{2}}\big]\frac{dx}{x} + \frac{1}{2}\int_{1}^{\infty}\big[x^{-\frac{s-1}{2}}-x^{\frac{-s}{2}}\big]\frac{dx}{x}$; since $\int_{1}^{\infty}x^{-a}\frac{dx}{x}=\frac{1}{a}$ for $a>0$, so $\frac{1}{2}\int_{1}^{\infty}\big[x^{-\frac{s-1}{2}}-x^{\frac{-s}{2}}\big]\frac{dx}{x} = \frac{1}{2}\bigg[\frac{1}{\frac{s-1}{2}}-\frac{1}{\frac{s}{2}}\bigg] = \frac{1}{s(s-1)}$.}, for $s>1$:

$$ \xi(s) = \Pi\bigg(\frac{s}{2}-1\bigg)\pi^{-\frac{s}{2}}\zeta(s) = \int_{1}^{\infty}\psi(x)\big[x^{\frac{s}{2}}+x^{\frac{1-s}{2}}\big]\frac{dx}{x} - \frac{1}{s(s-1)}.$$

\noindent Since $\psi(x)$ decreases more rapidly than any power of $x$ as $x \to \infty$, the integral in this formula converges for all $s$. This gives, therefore, another formula for $\zeta(s)$ which is valid for all $s$ other than $s=0,1$, i.e., it gives an alternate proof of the fact that $\zeta(s)$ can be analytically continued. Since both sides of the above equation are analytic, the equation holds for all $s$. Further, $\int_{1}^{\infty}\psi(x)\big[x^{\frac{s}{2}}+x^{\frac{1-s}{2}}\big]\frac{dx}{x} - \frac{1}{s(s-1)}$ is unchanged by the substitution $s = 1-s$, which proves the functional equation of the zeta function (Edwards, Section 1.7. "Second Proof of the Functional Equation" \cite{EHM}). We have the following series representation\footnote{Riemann states that this series representation of $\xi(s)$ as an even function of $s-\frac{1}{2}$ ``converges very rapidly".} for $\xi(s)$, for $s = \frac{1}{2}+bi$:
$$\xi(s) = \sum_{n=0}^{\infty}a_{2n}\bigg(s-\frac{1}{2}\bigg)^{2n}, \hphantom{...} \text{where } a_{2n} = 4\int_{1}^{\infty}\frac{d[x^{\frac{3}{2}}\psi'(x)]}{dx}x^{-\frac{1}{4}}\frac{(\frac{1}{2}\log x)^{2n}}{(2n)!}dx.$$

\newpage

\noindent Riemann applies Fourier\footnote{A ``master of evaluating and estimating definite integrals" (Edwards, Section 1.9), Riemann was also ``a master of Fourier analysis and his work in developing this theory must certainly be counted among his greatest contributions to mathematics" (Edwards, Section 1.12 ``Fourier Inversion" \cite{EHM}).} inversion to the formula (for Re $s>1$):
$$\frac{\log\zeta(s)}{s} = \int_{0}^{\infty} J(x) x^{-s-1}dx$$
$$\text{to obtain: }J(x) = \frac{1}{2\pi i}\int_{a-i\infty}^{a+i\infty}\log\zeta(s)x^s\frac{ds}{s}.$$
\noindent Here, it is understood that the improper integral above ``means the limit as $T\to\infty$ of the integral over the vertical line segment from $a-iT$ to $a+iT$" (Edwards, page 23 \cite{EHM}). The improper integral given above is only conditionally convergent and an ``order of summation" must be specified. Conditionally convergent integrals and series are very common in Fourier analysis, and it is always understood that such integrals and series are summed in their ``natural\footnote{Edwards provides some examples: $\sum_{n=-\infty}^{\infty}c_n e^{inx}$ means $\lim_{N\to\infty}\sum_{n=-N}^{N}c_n e^{inx}$, and\\ $\int_{-\infty}^{\infty}f(y)e^{iyx}dy$ means $\lim_{T\to\infty}\int_{-T}^{T}f(y)e^{iyx}dy$.} order" (Edwards, Section 1.12 ``Fourier Inversion" \cite{EHM}).

\noindent On page 86 of Riemann's \textit{Partielle Differentialgleichungen} (ed. Hattendorf) (1876) \cite{RBPD}, Riemann's use of ``Fourier's theorem" is revealed to be simply Fourier inversion. The reader is encouraged to review Paper XII of Riemann's Collected Papers, namely \textit{Abhandlungen der K\"{o}niglichen Gesellschaft der Wissenschaften zu G\"{o}ttingen, vol. 13}, ``On the representation of a function by a trigonometric series" \cite{RBCP}. A remark by Fourier threw new light on the topic of the analytic representation of arbitrary functions: for the trigonometric series $f(x)$ the coefficients $a_n$ and $b_n$ can be determined by the formula below and this method can also be applied if $f(x)$ is any arbitrary function:
$$ f(x) =
\begin{cases}
a_1 \sin x + a_2 \sin 2x + ...\\
+ \frac{1}{2}b_0 + b_1 \cos x + b_2 \cos 2x + ...
\end{cases}$$
$$a_n = \frac{1}{\pi}\int_{-\pi}^{\pi}f(x) \sin nx dx, \hphantom{...} b_n = \frac{1}{\pi}\int_{-\pi}^{\pi}f(x) \cos nx dx.$$

\noindent In Riemann's words, ``the applications of Fourier series are not restricted to research in the physical sciences. They are now also applied with success in an area of pure mathematics, number theory" (Riemann, Collected Papers, page 229 \cite{RBCP}).

\newpage
\noindent On page 5 \cite{RBU}, Riemann states that the number of roots of $\xi(t) = 0$ whose real parts lie between $0$ and $T$ is equal to the integral of $\frac{\xi'(s) ds}{2\pi i\xi(s)}$ around the boundary of the rectangle $\{0\leq \text{Re } s \leq 1, 0\leq \text{Im } s \leq 1,\}$ and this integral is equal to $\frac{T}{2\pi} \log\frac{T}{2\pi} - \frac{T}{2\pi}$ with a relative error of $\frac{1}{T}$. The fact that the vertical density of the roots $\rho$ is in some sense a constant times $\log\big(\frac{T}{2\pi}\big)$ was stated by Riemann without proof (Edwards, Section 1.10 ``The Product Representation of $\xi(s)$" \cite{EHM}). A proof of Riemann's estimate for the number of roots of $\xi(t) = 0$ whose real parts lie between $0$ and $T$ was first provided by von Mangoldt in \textit{Zur Verteilung der Nullstellen der Riemannschen Funktion $\xi(t)$} (1905) \cite{VMH}.

\noindent For Riemann, $\xi(t)$ is $\xi(s)$ with $s$ replaced by $\frac{1}{2}+ti$. On page 5 (transl. Wilkins, 1998), Riemann originally states: ``The number of roots of $\xi(t) = 0$, whose real parts lie between $0$ and $T$ is approximately: $= \frac{T}{2\pi} \log\frac{T}{2\pi} - \frac{T}{2\pi}$ because the integral $\int d \log \xi(t)$, taken in a positive sense around the region consisting of the values of $t$ whose imaginary parts lie between $\frac{1}{2}i$ and $-\frac{1}{2}i$ and whose real parts lie between $0$ and $T$, is (up to a fraction of the order of magnitude of the quantity $\frac{1}{T}$) equal to $\big(T \log\frac{T}{2\pi}-T\big)i$; this integral however is equal to the number of roots of $\xi(t) = 0$ lying within this region, multiplied by $2\pi i$. One now finds indeed approximately this number of real roots within these limits, and it is very probable that all roots are real" (Riemann \cite{RBU}). As Bombieri remarks in \textit{Problems of the Millennium: The Riemann Hypothesis} \cite{BE}, the Riemann Hypothesis is the statement:
\begin{center}
\textit{All zeros of the function $\xi(t)$ are real}.
\end{center}
\vspace{10 pt}
\begin{center}
    The Riemann Hypothesis remains unproven to this day.
\end{center}

\newpage
\thispagestyle{empty}
\section{From the Prime Number Theorem to the M\"{o}bius Reformulation of RH}

\noindent The conjecture of Gauss and Legendre\footnote{In ``Historical Introduction", \textit{Introduction to Analytic Number Theory} \cite{ATM}, Apostol traces the history of $\pi(n)$ back to Legendre (in \textit{Essai sur la Theorie des Nombres}, 1798 \cite{LAM}) and Gauss (in ``Letter to Encke", dated 24, December, 1849 \cite{GL}). Apostol provides a table with columns $x$, $\pi(x)$, $\frac{x}{\log x}$, and $\frac{\pi(x)}{\frac{x}{\log x}}$, suggesting that though the distribution of the prime numbers is irregular, ``by examining large blocks of primes one finds that their average distribution seems to be quite regular" (page 8, Apostol \cite{ATM}). Thus, Gauss and Legendre correctly conjecture that as $x \to \infty$, $\frac{\pi(x)}{\frac{x}{\log x}} \to 1$, but lack the means to provide a proof. In 1851, Chebyshev \cite{CH} proved that if the ratio does tend to a limit, then this limit must be $1$.} remained unproven till in 1896, Hadamard (in \textit{Sur la distribution des z\'{e}ros de la fonction $\zeta(s)$ et ses cons\'{e}quences arithm\'{e}tiques} \cite{HJS}) and C. J. de la Vall\'{e}e Poussin (in \textit{Recherches analytiques sur la th\'{e}orie des nombres premiers} \cite{VPCJ}), using complex analytic methods, independently proved what is now referred to as the ``Prime Number Theorem" (PNT), i.e.,
$$ \lim_{x \to \infty} \frac{\pi(x)}{\frac{x}{\log x}} = 1$$

\noindent In 1949, Selberg (\textit{An elementary proof of the prime number theorem} \cite{SA}) and Erd\"{o}s (\textit{On a new method in elementary number theory which leads to an elementary proof of the prime number theorem} \cite{EP}) discovered ``elementary" proofs of the prime number theorem, without using $\zeta (s)$ or complex function theory.

\noindent The idea of furnishing an elementary proof of the prime number theorem ``caused a sensation in the mathematical world" (Apostol, page 9 \cite{ATM}). Before Selberg and Erd\"{o}s published their findings, Hardy addressed the Mathematical Society of Copenhagen in 1921, ruling out the possibility of an elementary proof of the prime number theorem, asserting that ``some theorems ``lie deep" and others nearer to the surface", and that ``if anyone produces an elementary proof of the prime number theorem, he will show that these views are wrong, that the subject does not hang together in the way we have supposed, and that it is time for the books to be cast aside and for the theory to be rewritten" (Hardy, as quoted by Goldfield in \textit{The Elementary Proof of the Prime Number Theorem: An Historical Perspective}, 2003 \cite{GD}).

\newpage

\noindent In Landau's \textit{Handbuch der Lehre von der Verteilung der Primzahlen} \cite{LE}, or simply \textit{``Handbuch"}, hailed by Hardy (in the obituary of Landau he wrote for the London Mathematical Society \cite{HGH}) as the first time the analytic theory of numbers is presented ``not as a collection of a few beautiful scattered theorems, but as a systematic science", one may flip through the first few pages to see that Landau begins this systematic science with two simple functions, namely the M\"{o}bius function, $\mu$, and the Liouville function, $\lambda$.

\noindent The M\"{o}bius function, $\mu(n)$, is defined as $(-1)^k$, if $n$ is the product of $k$ \textit{distinct} primes and $0$ otherwise; the Liouville Function $\lambda (n)$, on the other hand, is defined as $(-1)^k$, if $n$ is a product of $k$ primes (including multiplicities). Remarkably, $\mu$ and $\lambda$ are closely related to $\zeta$. On page 12 of \textit{The Theory of the Riemann Zeta-Function} \cite{TEC}, Titchmarsh provides the following relationships between $\zeta$ and $\mu$ (equation 1.1.4, Titchmarsh \cite{TEC}), and $\zeta$ and $\lambda$ (equation 1.2.11, Titchmarsh \cite{TEC}), for $s=a+ib$, with $a>1$, as follows, the second identity attributed to Lehman (\textit{On Liouville's Function}, 1960 \cite{LRS}):
$$\frac{1}{\zeta(s)} = \sum_{n=1}^{\infty}\frac{\mu(n)}{n^s} \quad \quad \quad \quad \frac{\zeta(2s)}{\zeta(s)} = \sum_{n=1}^{\infty} \frac{\lambda (n)}{n^s} $$
\noindent In Chapter 3 ``Averages of arithmetical functions", \textit{Introduction to Analytic Number Theory} \cite{ATM}, Apostol says that the statement ``$\mu(n)$ has average order $0$", i.e., ``$\sum_{n \to \infty}^{}\frac{\mu(n)}{n}$ converges and has sum $0$", is equivalent to the PNT. Further, as referenced by Spencer-Brown in Appendix 9 \cite{SBG}, by Edwards in Section 12.3, ``Miscellany" \cite{EHM}, by Apostol in Exercise 4 of Chapter 13 (page 301) \cite{ATM}, and by Derbyshire as Theorem 15-2, ``Big Oh\footnote{We say $f(x)$ is $O(g(x))$ or "big-O of $g(x)$" if there are positive constants $C$ and $x_0$ such that $|f(x)| \leq C|g(x)|$ whenever $x>x_0$.} and M\"{o}bius Mu" \cite{DJ}, RH is equivalent to $M(x) = O(x^{\frac{1}{2}+\epsilon}), \hphantom{...} \forall \epsilon > 0$, where $M(x) = \sum_{n \leq x}^{} \mu(n)$.

\noindent Similarly, Theorem 1.4 in \textit{The Riemann Hypothesis} (Borwein et al, 2006) \cite{BP}, states that the PNT is equivalent to the statement $\lim_{n \to \infty}\frac{\lambda(1) + \lambda(2) + ... + \lambda(n)}{n} = 0$. Further, as also referenced by Spencer-Brown on page 209 \cite{SBG}, and Borwein et al on page 6 (Theorem 1.2) \cite{BP}, RH\footnote{Spencer-Brown reiterates here that upon subsituting $n$ for $n^{\frac{1}{2}+\epsilon}$ one obtains the much weaker PNT.} is equivalent to $\lim_{n \to \infty}\frac{\lambda(1) + \lambda(2) + ... + \lambda(n)}{n^{\frac{1}{2}+\epsilon}} = 0$, for any fixed $\epsilon > 0$. Both of these results are attributed to Landau in his dissertation (1899) \cite{LED}.

\newpage

\noindent Let $dM$ be the Stieltjes measure such that if $\frac{1}{\zeta(s)} = \sum_{n=1}^{\infty}\frac{\mu(n)}{n^s}, (\text{for Re } s>1)$,
$$ \frac{1}{\zeta(s)} = \int_{0}^{\infty} x^{-s} d M(x),$$
\noindent where $M(x) = \int_{0}^{x} dM$ is a step function which is zero at $x=0$, which is constant except at positive integers, and which has a jump\footnote{The value of $M$ at a jump is by definition $\frac{1}{2}\big[M(n-\epsilon)+M(n+\epsilon)\big] = \sum_{j=1}^{n-1} \mu(j) + \frac{1}{2} \mu(n)$.} of $\mu(n)$ at $n$ (Edwards, Section 12.1, ``The Riemann Hypothesis and the Growth of $M(x)$" \cite{EHM}). Integration by parts for Re $s>1$, gives:
$$\frac{1}{\zeta(s)} = $$
$$\int_{0}^{\infty} d\big[x^{-s} M(x)\big] - \int_{0}^{\infty} M(x) d(x^{-s}) = \lim_{X \to \infty}\bigg[X^{-s} M(X) + s \int_{0}^{X} M(x) x^{-s-1} dx\bigg]$$
$$ = s \int_{0}^{\infty} M(x) x^{-s-1} dx$$

\noindent because the inequality $|M(x)| \leq x$ implies that $x^{-s} M(x) \to 0$ as $x \to \infty$ and $\int_{0}^{\infty} M(x) x^{-s-1} dx$ converges, both provided Re $s>1$. If $M(x)$ grows less rapidly than $x^a$ for some $a>0$, then this integral for $\frac{1}{\zeta(s)}$ converges for all $s$ in the halfplane $\{ \text{Re } (a-s) < 0\} = \{\text{Re } s>a\}$. By analytic continuation, the function $\frac{1}{\zeta(s)}$ is analytic in this halfplane. Since $\frac{1}{\zeta(s)}$ has poles on the line Re $s = \frac{1}{2}$,
\begin{center}
$M(x)$ \textit{does not grow less rapidly than $x^a$ for any $a<\frac{1}{2}$}.  
\end{center}

\noindent Moreover, it shows that \textit{in order to prove the Riemann hypothesis, it would suffice to prove that $M(x)$ grows less rapidly than $x^{\frac{1}{2}+\epsilon}$ for all $\epsilon > 0$} (Edwards, Section 12.1, page 260-261 \cite{EHM}). Littlewood in \textit{Quelques cons\'{e}quences de l'hypoth\`{e}se que la fonction $\zeta(s)$ n'a pas de z\'{e}ros dans le demi-plan Re $(s) > \frac{1}{2}$} (1912) \cite{LJE} proved that this sufficient condition for the Riemann hypothesis is also necessary, hence proving the following \textbf{theorem}, attributed to Hardy:

\begin{center}
    The Riemann Hypothesis (RH) is equivalent to the statement that
\end{center}
 $$\text{for every } \epsilon > 0, \frac{M(x)}{x^{\frac{1}{2}+\epsilon}} \text{ approaches zero as  }x \to \infty.$$

\noindent This theorem will be very important for the rest of the paper. On the one hand, it is the source of speculation for RH and coin-tossing and on the other hand, it is the central aspect with which Spencer-Brown writes Appendix 9 of Laws Of Form.

\newpage
\thispagestyle{empty}
\section{Appendix 9 of ``Laws Of Form" and Denjoy's Probabilistic Interpretation \cite{EHM}}
\noindent An argument for the plausibility of RH can be found in Chapter 12, ``Miscellany" of Edwards' \textit{Riemann's Zeta Function}, in Section 12.3, ``Denjoy's Probabilistic Interpretation\footnote{In pages 195 - 196 of \textit{Probabilit\'{e}s confirmant l'hypoth\`{e}se de Riemann sur les z\'{e}ros de $\zeta(s)$} \cite{DA}, we find Denjoy's original statement: ``En outre $\mu(1) = 1$.
$$\text{Posons}, \hphantom{...} \Delta(n) = \sum_{1}^{n} \mu(i). g(s) = \sum_{1}^{\infty}\Delta(n)[n^{-s}-(n+1)^{-s}] \sim s \sum \Delta(n) n^{-1-s},$$
$$ \mid s \mid^{-1} \mid g(s) \mid \leq \sum \mid \Delta (n) \mid n^{-1-\sigma}$$
Il suffrait donc de prouver que $\overline{\lim_{n\to\infty}} \log \Delta(n)/\log n  = \frac{1}{2}$ (l'in\'{e}galit\'{e} $< \frac{1}{2}$ est impossible) pour que l'hypoth\`{e}se de Riemann soit justifi\'{e}e." Also, ``Si $s = \sigma + it$, $\mid \sigma$ et $t$ \'{e}tant r\'{e}els, la s\'{e}rie et le produit infini convergent absolument pour $\sigma > 1$" (Denjoy, page 195 \cite{DA}).

} of the Riemann Hypothesis" \cite{EHM}. On page 268, Edwards reasons as follows: ``Suppose an unbiased coin is flipped a large number of times, say $N$ times. By the de Moivre-Laplace theorem the probability that the number of heads deviates by less than $KN^{\frac{1}{2}}$ from the expected number of $\frac{1}{2}N$ is nearly equal to $\int_{-(\frac{2K^2}{\pi})^{\frac{1}{2}}}^{(\frac{2K^2}{\pi})^{\frac{1}{2}}} \exp (-\pi x^2) dx$ in the sense that the limit of these probabilities as $N \to \infty$ is equal to this integral. Thus if the total number of heads is subtracted from the total number of tails, the probability that the resulting number is less than $2KN^{\frac{1}{2}}$ in absolute values is nearly equal to $2 \int_{0}^{\frac{2K^2}{\pi}^{\frac{1}{2}}} \exp (-\pi x^2) dx$. The fact that this approaches $1$ as $N \to \infty$ can be regarded as saying that \textit{with probability one the number of heads minus the number of tails grows less rapidly than $N^{\frac{1}{2}+\epsilon}$}" (Edwards, page 268, Section 12.3 \cite{EHM}).

\noindent If we consider a very large square-free integer $n$, then $\mu(n) = \pm 1$. ``But then the evaluation of $M(x)$ would be like flipping a coin once for each square-free integer less than $x$ and subtracting the number of heads from the number of tails. It was shown above that for any given $\epsilon > 0$ the outcome of this experiment for large\footnote{``The number of flips goes to infinity as $x \to \infty$ because, among other reasons, there are infinitely many primes, hence \textit{a fortiori} infinitely many square-free integers (products of distinct primes)" (original footnote, Edwards, page 268).} number of flips is, with probability nearly one, less than the number of flips raised to the power $\frac{1}{2} + \epsilon$ and \textit{a fortiori} less than $x^{\frac{1}{2}+\epsilon}$" (Edwards, page 268, Section 12.3 \cite{EHM}).

\newpage

\noindent On page 268, Section 12.3, ``Denjoy's Probabilistic Interpretation of the Riemann Hypothesis" \cite{EHM}, Edwards, however, goes on to write, ``it is perhaps plausible to say that successive evaluations of $\mu (n) = \pm 1$ are ``independent" since knowing the value of $\mu (n)$ for one $n$ would not seem to give any\footnote{An exception to this statement is that for any prime $p$, $\mu (pn)$ is either $- \mu (n)$ or zero. However, this principle can only be applied once for any $p$ because $\mu(p^2 n) = 0$ and this ``information" really says little more than that $\mu$ is determined by a formula and is not, in fact, a random phenomenon (original footnote by Edwards).} information about its values for other values of $n$."

\noindent As pointed out by Spencer-Brown, it would be remiss to consider successive evaluations of $\mu(n) = \pm 1$ as independent of one another. In \textit{Essai sur la th\'{e}orie des nombres} \cite{LAM}, Legendre reformulates the prime counting function $\pi (n)$ as a function of $\pi (\sqrt{n})$ using the M\"{o}bius function $\mu$. Spencer-Brown rectifies this formula for $d = 1$ and offers a simple recursive formula to compute $\mu(n)$, thus showing that subsequent values of $\mu(n) = \pm 1$ are \textbf{not} ``independent".

\noindent The Spencer-Brown Formula\footnote{``There are two Spencer-Brown formulae $Mu (n)$. They are both algebraic and connected by the sign of equality$: `='$. It should be observed that both formulae are sensitive to parentheses. If we remove the parentheses from the first formula, the minus sign distributes and converts the formula into the second. Conversely, if we add parentheses and a minus sign the second formula converts once more into the first. The formula was designed specifically to prove RH. It is a modification of procedures employed by Legendre and Sierpi\'{n}ski" (J. M. Flagg in \textit{Conversations with GSB, Oct 8, 2012}).} M\"{o}bius $\mu(n)$ is presented as follows: for $n>1$,

$$\mu (n) = - \Big( n-1 + \sum_{d=2}^{n-1} \mu (d)\left[ \frac{n}{d} \right] \Big) = 1 - n - \sum_{d=2}^{n-1} \mu (d)\left[ \frac{n}{d} \right],$$
\begin{center}
where the bracketed variable $\big[\frac{n}{d}\big]$ denotes the integer part of this fraction.
\end{center}

\newpage

\noindent The sixth English edition of ``Laws of Form" (2015) \cite{SBG} includes George Spencer-Brown's most recent version of ``Appendix 9: A proof of Riemann's hypothesis using Denjoy's equivalent theorem". It begins with a ``Description of the proof", which states that ``The Riemann hypothesis is true if and only if the numbers of positive and negative signs of $\mu (n)$ are asymptotically\footnote{Historically credited to Paul Bachmann in \textit{Die Analytische Zahlentheorie} ``Analytic Number Theory", 1894, the big $O$ notation currently used to describe orders of magnitude was first clearly defined by Landau. "More specifically, the $f(x) \sim g(x)$ relation for two functions $f, g : \mathbb{R} \to \mathbb{R}$ means that $\lim_{x\to\infty}\big( \frac{f(x)}{g(x)}\big) = 1$.  Under these conditions, the two functions are said to be \textit{asymptotically} equivalent (or simply \textit{asymptotic}). Alternatively, we say that $f(x)$ is $O(g(x))$ (pronounced "big-O of $g(x)$") if there are positive constants $C$ and $x_0$ such that $|f(x)| \leq C|g(x)|$ whenever $x>x_0$. It is not too difficult to see that $f(x) \sim g(x)$ implies that $f(x)$ is $O(g(x))$ and $g(x)$ is $O(f(x))$" (Tou, in \textit{Math Origins: Orders of Growth} \cite{TE}).} equal", referencing Denjoy's probabilistic interpretation\footnote{\textit{Probabilit\'{e}s confirmant l'hypoth\`{e}se de Riemann sur les z\'{e}ros de $\zeta(s)$. Note $(*)$ de M. Arnaud Denjoy, Membre de l'Acad\'{e}mie, in ``Conclusions tir\'{e}es d'une forme donn\'{e}e au d\'{e}veloppement de $\frac{\zeta(2s)}{\zeta(s)}$}" \cite{DA}}, that ``RH is equivalent to the proposition that any square-free number, taken at random, has an equal probability of containing an odd or an even number of (different) prime divisors" (page 203, ``Laws of Form", 2015 \cite{SBG}). We note that Spencer-Brown is undoubtedly referring to the full technical statement that $M(x) = O(x^{\frac{1}{2}+\epsilon})$.

\noindent The Legendre/Spencer-Brown formula for the prime counting function $\pi (n)$ is a recursive function that calculates the number of primes $\leq n$ by referring to $\pi (\sqrt{n})$. The formula only uses primes up to $\sqrt{n}$, what Spencer-Brown calls ``small" primes, to sieve the remaining primes up to $n$, what Spencer-Brown calls the ``large" primes. The Legendre formula, rectified by Spencer-Brown, is presented as follows:

$$ \pi (n) = \pi ( \sqrt{n}) + \sum_{d = 1}^{\sqrt{n}} \mu (d) \left[ \frac{n}{d} \right],$$

\noindent Specifically for $d = 1$, Spencer-Brown rectifies $\left[\frac{n}{d}\right]$ to $\left[\frac{n-1}{d}\right]$. Here, the bracketed variable $\left[ \frac{n}{d} \right]$ denotes the integer part of $\frac{n}{d}$. In general, $d$ consists of primes up to $\sqrt{n}$, and all possible products of these primes. Since $\mu(d)=0$ when there are repeated prime factors, we only consider those $d$ that factor into single primes up to $\sqrt{n}$.

\newpage

for $n = 100$,
$$ \pi (100) = \pi (10) + \sum_{d = 1}^{\sqrt{100}} \mu (d) \left[ \frac{100}{d} \right]. $$

\noindent The small primes are all primes up to $\sqrt{100} = 10$, which are $2, 3, 5,$ and $7$. They combine to form products that divide $100$. They can be organized by collecting the doubles, triples, and so on:
$$\overbrace{2, 3, 5, 7}, \overbrace{(2.3), (2.5), (2.7), (3.5), (3.7), (5.7)}, \overbrace{(2.3.5), (2.3.7), (2.5.7), (3.5.7)}, (2.3.5.7)$$

$$\text{Hence, }\sum_{d = 1}^{\sqrt{100}} \mu (d) \left[ \frac{100}{d} \right]$$
$$ = \Bigg( \mu (1) \left[ \frac{100-1}{1} \right] + \mu (2) \left[ \frac{100}{2} \right] + \mu (3) \left[ \frac{100}{3} \right] + \mu (5) \left[ \frac{100}{5} \right] + \mu (7) \left[ \frac{100}{7} \right]$$
$$+ \mu (2.3)\left[\frac{100}{2.3}\right] + \mu (2.5)\left[\frac{100}{2.5}\right] + \mu(2.7)\left[\frac{100}{2.7}\right] + \mu(3.5)\left[\frac{100}{3.5}\right] + \mu(3.7)\left[\frac{100}{3.7}\right]$$
$$+ \mu(5.7)\left[\frac{100}{5.7}\right] + \mu (2.3.5) \left[ \frac{100}{2.3.5} \right] + \mu (2.3.7) \left[ \frac{100}{2.3.7} \right] + \mu (2.5.7) \left[ \frac{100}{2.5.7} \right]$$
$$+ \mu (3.5.7) \left[ \frac{100}{3.5.7} \right] + \mu (2.3.5.7) \left[ \frac{100}{2.3.5.7} \right] \Bigg)$$

$$ = \Bigg( \mu (1) \left[ \frac{99}{1} \right] + \mu (2) \left[ \frac{100}{2} \right] + \mu (3) \left[ \frac{100}{3} \right] + \mu (5) \left[ \frac{100}{5} \right] + \mu (7) \left[ \frac{100}{7} \right] $$
$$ + \mu (2.3) \left[ \frac{100}{6} \right] + \mu (2.5) \left[ \frac{100}{10} \right] + \mu (2.7) \left[ \frac{100}{14} \right] + \mu (3.5) \left[ \frac{100}{15} \right] + \mu (3.7) \left[ \frac{100}{21} \right]$$
$$ + \mu (5.7) \left[ \frac{100}{35} \right] + \mu (2.3.5) \left[ \frac{100}{30} \right] + \mu (2.3.7) \left[ \frac{100}{42} \right] + \mu (2.5.7) \left[ \frac{100}{70} \right]$$
$$+ \mu (3.5.7) \left[ \frac{100}{105} \right] + \mu (2.3.5.7) \left[ \frac{100}{210} \right] \Bigg)$$

$$ = \Bigg( +\left[ \frac{99}{1} \right] -\left[ \frac{100}{2} \right] -\left[ \frac{100}{3} \right] -\left[ \frac{100}{5} \right] -\left[ \frac{100}{7} \right] $$
$$ +\left[ \frac{100}{6} \right] +\left[ \frac{100}{10} \right] +\left[ \frac{100}{14} \right] +\left[ \frac{100}{15} \right] +\left[ \frac{100}{21} \right] +\left[ \frac{100}{35} \right] $$
$$ -\left[ \frac{100}{30} \right] -\left[ \frac{100}{42} \right] -\left[ \frac{100}{70} \right] -\left[ \frac{100}{105} \right] +\left[ \frac{100}{210} \right] \Bigg)$$
$$ = + 99 - 50 - 33 - 20 - 14$$
$$ + 16 + 10 + 7 + 6 + 4 + 2$$
$$ - 3 - 2 - 1 - 0 + 0$$
$$ = 21$$
\noindent Hence, according to LSB (Legendre/Spencer-Brown), $\pi (100) = \pi (10) + 21$. But we know $\pi(10) = 4$, when we listed out all the small primes $2, 3, 5$, and $7$.
$$ \pi(100) = \pi(10) + 21$$
$$ = 4 + 21 = \textbf{25}.$$

\noindent The Legendre/Spencer-Brown formula can be tested for any $n$ and the formula ``works correctly because its summation term yields the number of numbers $\leq n$ that are not struck out by the Eratosthenes procedure of striking out those of them that are divisible by a prime $q$ that is `small' in relation to $n$, i.e., such that $2 \leq q \leq n^{\frac{1}{2}}$" (page 203, Spencer-Brown, 2015 \cite{SBG}). In every application of the Legendre procedure, the small primes, i.e., the primes up to $\sqrt{n}$, are needed to count the large primes, i.e., primes up to $n$. Two important observations:

\noindent When the divisor $d$ exceeds $n$, the `floored' quotient $\left[\frac{n}{d}\right]$ returns $0$ and so $f(d) = \mu (d) \left[\frac{n}{d}\right] = \mu (d).0 = 0$. We saw this with divisors $3.5.7 = 105$ and $2.3.5.7 = 210$ in the Legendre/Spencer-Brown demonstration for $n = 100$ on the previous page.

\noindent When a square-free divisor $d$ exceeds $\frac{n}{2}$ but remains $\leq n$, however, the quotient $\left[\frac{n}{d}\right]$ yields $1$, and $f(d) = \mu (d) \left[\frac{n}{d}\right]$ is simply $\mu (d)$. This important observation leads Spencer-Brown to call the square-free numbers $\leq \frac{n}{2}$ the ``upper section", and the square-free numbers $> \frac{n}{2}$ but $\leq n$ the ``lower section". On page 205, Spencer-Brown writes, ``in the lower section the sum of the terms is exactly the sum of the $\mu(d)$ in the section. Call a series of such terms an LSB series" (page 205, Appendix 9 \cite{SBG}).

\newpage

\noindent ``The fact that with $d$ unrestricted, the formula $\sum \mu(d)\big[\frac{n}{d}\big]= 1$ is true for all $n$ was first noted by Meissel in \textit{Observationem quaedam in theoria numerorum}, Berlin 1850, and proved by Sierpi\'{n}ski in \textit{Elementary theory of numbers}, Warsaw 1964, pp 180 and 181" (Spencer-Brown, footnote on page 204 \cite{SBG}). The point here is that instead of summing the small primes and their products as in the Legendre summation, one can sum all $d$ from $1$ to $n$ so one has $\sum_{d=1}^{n} \mu(d)\big[\frac{n}{d}\big] = 0$, by using Spencer-Brown’s adjustment. On page 204, he provides an illustration of the use of Legendre's formula to calculate $\pi(n)$ for $n=20$:

\begin{multicols}{2}
\noindent $d$ \quad \quad $f(d) = \mu(d)\big[\frac{n}{d}\big]$\\
$1^*$ \quad \quad $+20$\\
$2^*$ \quad \quad $-10$\\
$3^*$ \quad \quad $-6$\\
$5$\hphantom{.} \quad \quad $-4$\\
$6^*$ \quad \quad $+3$\\
$7$\hphantom{.} \quad \quad $-2$\\
$10$ \quad \quad $+2$\\
$11$ \quad \quad $-1$\\
$13$ \quad \quad $-1$\\
$14$ \quad \quad $+1$\\
$15$ \quad \quad $+1$\\
$17$ \quad \quad $-1$\\
$19$ \quad \quad \underline{$-1$}\\
\hphantom{.} \quad \quad $\sum +1$\\
%\vfill\null
\columnbreak

\noindent $\sum \mu(d^*)\big[\frac{n}{d^*}\big] = 7$\\
$7-1+\pi(n^{\frac{1}{2}}) = 8 = \pi(20)$

\noindent The $(d^*)$ are the denominators with no large prime in their decomposition. The small primes $2, 3$ must be known explicitly, then the number of large primes $5, 7, 11, 13, 17, 19$ can be calculated without any of them being identified.
\end{multicols}

\noindent In Spencer-Brown's words: ``Since $1$ is not struck out by the sieve of Eratosthenes and is also included in the count of ``primes" calculated by the section $\sum \mu(d) \big[\frac{n}{d}\big]$, the count must be reduced by one in either case, and then to get the complete answer the number of small primes $(q)$ used as strikers must be added to the total" (Spencer-Brown, on page 204). ``Notice I have used all $(d)$ that yield an $f(d)$ other than zero, and the sum of these, for any $n$, must always be $1$, since only one number, $1$ itself, remains unstruck if we use all the primes" (Spencer-Brown, on page 204 \cite{SBG}).

\newpage

\noindent ``We thus further rectify the procedure by making the following change: use $f(d) = \mu(d) \big[\frac{n-1}{d}\big]$ for $d = 1$, and use $f(d) = \mu(d) \big[\frac{n}{d}\big]$ for all other values of $d$" (Spencer-Brown, page 205 \cite{SBG}). This rectification ensures:
$$\sum_{d=1}^{n} f(d) = 0, \hphantom{... ...} \text{where } f(d) = \mu(d)\floor*{\frac{n}{d}}, \hphantom{..}\text{and for }d = 1, f(d) = \mu(d)\floor*{\frac{n-1}{d}}.$$

\noindent ``My Denjoy proof, in summary, runs as follows:\\
1. What Professor Denjoy showed is that the RH is equivalent to the proposition that the number of primes in a square-free $d \leq n$ of any size is even or odd with equal probability.\\
2. I rectify Legendre's method of counting large primes in $n$ and then corrupt it to give the answer zero for all $(n)$.\\
3. I split the rectified Legendre terms into two sections, upper and lower, so that the lower sections eventually include all values of $\mu(d)$ for square-free $(d)>1$.\\
4. I show that the average algebraic sum, i.e., the sum divided by the number of LSB terms displayed, in each section varies around and is asymptotic to zero as $n$ for the $f(d)$ terms increases without limit.\\
5. Since the lower sections eventually include the values of $\mu(d)$ for all square-free $(d)>1$, and their signs are by the previous proposition equiprobable in the limit, the RH, quod erat demonstrandum, must be true.\\
6. In addition, since the upper sections eventually contain all the values of $\mu(d)$ for square-free $(d)$ but magnified by various factors ranging from $2$ to $n-1$, that are independent of the signs of $\mu(d)$, and the average differences between the plus and minus values of these magnified terms also tend to zero as $n$ increases, this fact constitutes a second proof of Riemann's hypothesis, since if an average set of magnified differences tends to zero, then the average of the same set of differences unmagnified must also tend to zero" (page 214, Appendix 9 \cite{SBG}).

\noindent The negative feedback property of $\mu(n)$ is used by Spencer-Brown to argue that the cumulative M\"{o}bius function varies around zero and that its asymptotics are sufficient, by Denjoy's interpretation, to prove the Riemann Hypothesis (RH), which is equivalent to the relation $M(x) = O(x^{\frac{1}{2}+\epsilon}), \forall \epsilon>0$, where $M(x) = \sum_{n \leq x}^{} \mu(n)$.

%\noindent The approach of Spencer-Brown deserves further investigation. We have, in this paper, placed his approach in the context of the evolution of methods for understanding RH, and in particular in relation to the reformulation of RH in terms of the behaviour of the cumulative Mobius Function $M(n)$. In this light, we know that to prove RH one needs to show that $M(n)$ is of the order $O(\sqrt{n} n^{\epsilon})$. On the other hand, Spencer-Brown has given elegant formulas for finding the M\"{o}bius function in terms of its previous values, and has used the Meissel-Legendre formula to show the behaviour of $M(n)$ in terms of feedback and recursion. Spencer-Brown points out that the recursion implicit in the Spencer-Brown/Meissel-Legendre formula makes $M(n)$ self-regulating and that it is this self-regulation that drives its asymptotic behaviour. We want to know if these methods do lead to a formal proof of the needed behaviour of $M(n)$. Spencer-Brown has said that this is so.

\newpage

\noindent We have, in this paper, placed Spencer-Brown's approach to the Riemann Hypothesis in the context of the evolution of the problem, and in particular in relation to the reformulation of RH in terms of the behaviour of the cumulative M\"{o}bius Function $M(n)$. In this light, we know that to prove RH is to show that $M(n)$ is of the order $O(\sqrt{n} n^{\epsilon})$. On the other hand, Spencer-Brown has given a deep formula for finding the M\"{o}bius function in terms of its previous values. This formula is
$$ \mu(n) =  - \Bigg(\sum_{d=1}^{n-1} \mu(d)\floor*{\frac{n}{d}}\Bigg),$$
where for $d=1$, we make the Spencer-Brown adjustment, rectifying $\big[\frac{n}{d}\big]$ to $\big[\frac{n-1}{d}\big]$.

\noindent The key point about this formula is that it is a negative feedback for the next value of the M\"{o}bius function with respect to all of its previous values. The formula shows that the M\"{o}bius function is far from being random and indeed that any run of strictly positive values or any run of strictly negative values will be turned around by the negative feedback. This means that a cumulative sum of the values of the M\"{o}bius function will exhibit highly controlled oscillation about zero, with widening swings as the value of $n$ in $M(n)$ grows. The recursion implicit in the Spencer-Brown/Legendre formula makes $M(n)$ self-regulating and it is this self-regulation that drives its asymptotic behaviour. Spencer-Brown's fundamental observation is that the M\"{o}bius function $\mu (n)$, while not random, behaves in a more constrained way than a random series due to the feedback inherent in its dependence on its previous values. To quote Spencer-Brown,  ``$\cdots$ the successive nonzero values of $\mu(d)$ behave like the falls of a magic coin that, from the moment it is struck, remembers exactly how many times it has fallen with one side or the other uppermost, and whenever one side exceeds the other, biases itself towards the other until the excess is eliminated" (page 211, Appendix 9 \cite{SBG}). This means that the asymptotics for a random series still apply to the M\"{o}bius function and so the Riemann Hypothesis (RH) follows.

\noindent What should be observed is that the very technology of Laws of Form, the mark of distinction, is itself the perfect medium to employ for the Riemann Hypothesis, as every distinction cleaves its domain in two or what is the same, $1/2$, or two halves. Such spaces are always in balance and objects within them must behave according to a maturity of chances hypothesis. This is due to feedback and feed forward everywhere present as shown in the cascade. The cascade is the simplest possible arithmetic one can develop that exhibits this behaviour and manifests that behaviour throughout all aspects of the Riemann Hypothesis. \\

%\newcounter{pnumb}                     	% Save page number of intro material
%\setcounter{pnumb}{\value{page}}       	% so that numbering will be consecutive.
%\mainmatter                            	% Begin body
%\setcounter{page}{\value{pnumb}}       	% with consecutive page numbering
\backmatter
\thispagestyle{empty}
\section*{Appendix \texorpdfstring{$\pm1$}{} -- The Spencer-Brown M\"{o}bius Cascade}
\noindent We are presented with The Spencer-Brown Formula\footnote{$$( \text{for } n>1 ),\hphantom{...} \mu (n) = - \Big( n-1 + \sum_{d=2}^{n-1} \mu (d)\left[ \frac{n}{d} \right] \Big) = 1 - n - \sum_{d=2}^{n-1} \mu (d)\left[ \frac{n}{d} \right],$$
where the bracketed variable $\left[ \frac{n}{d} \right]$ denotes the integer part of $\frac{n}{d}$.} M\"{o}bius $\mu(n)$, which can be seen to follow as a consequence of the Legendre/Spencer-Brown formula for $\pi(n)$. The Spencer-Brown Formula M\"{o}bius $\mu(n)$ (J. M. Flagg, \textit{Conversations with GSB, Oct 8, 2012)}:

$$\mu(n) = - \Big( n-1 + \sum_{d=2}^{n-1} \mu (d) \floor*{\frac{n}{d}} \Big) = 1 - n - \sum_{d=2}^{n-1} \mu (d)\floor*{\frac{n}{d}}$$\\
\medskip
$$ =  - \Bigg(\sum_{d=1}^{n-1} \mu(d)\floor*{\frac{n}{d}}\Bigg),$$
\noindent since for $d=1$, we make the Spencer-Brown adjustment, rectifying $\big[\frac{n}{d}\big]$ to $\big[\frac{n-1}{d}\big]$.

\noindent Start with $\mu(1) = \mathbf{+1}$.\\
\bigskip
\noindent $\mu(2) = - \big(\mu(1)\big[\frac{2-1}{1}\big]\big) = \mathbf{-1}$.\\
\bigskip
\noindent $\mu(3) = - \big(\mu(1)\big[\frac{3-1}{1}\big] + \mu(2)\big[\frac{3}{2}\big]\big) = - (2 - 1) = \mathbf{-1}$.\\
\bigskip
\noindent $\mu(4) = - \big(\mu(1)\big[\frac{4-1}{1}\big] + \mu(2)\big[\frac{4}{2}\big] + \mu(3)\big[\frac{4}{3}\big]\big) = - (3 - 2 - 1) = \mathbf{0}$.\\
\bigskip
\noindent $\mu(5) = - \big(\mu(1)\big[\frac{5-1}{1}\big] + \mu(2)\big[\frac{5}{2}\big] + \mu(3)\big[\frac{5}{3}\big] + \mu(4)\big[\frac{5}{4}\big] \big) = - (4 - 2 - 1 + 0) = \mathbf{-1}$.\\
\bigskip
\noindent $\mu(6) = - \big(\mu(1)\big[\frac{6-1}{1}\big] + \mu(2)\big[\frac{6}{2}\big] + \mu(3)\big[\frac{6}{3}\big] + \mu(4)\big[\frac{6}{4}\big] + \mu(5)\big[\frac{6}{5}\big] \big) = - (5 - 3 - 2 + 0 - 1) = \mathbf{+1}$.

\newpage

\noindent Spencer-Brown's cascade introduced in Appendix 9 \cite{SBG} is a tabular arrangement:  
The first column lists all numbers $d\leq n$. The second column lists $f(d) = \mu(d)\floor*{\frac{n}{d}}$ (where, for $d=1, \frac{n}{d}$ is rectified to $\frac{n-1}{d}$).
The third column lists $\sum f(d)$, what Spencer-Brown calls the ``running total", and whose final term is always $0$.

\noindent for $n = 2$,\\
$\hphantom{..}d \quad \hphantom{...} f(d) \hphantom{..} \sum f(d) \hphantom{...}$\\ 
$\hphantom{..}1 \quad \hphantom{...}  + 1 \quad + 1 $\\
$\mu (2) = - 1 \quad \quad0$

\noindent for $n=3$,\\
$\hphantom{..}d \quad \hphantom{...} f(d) \hphantom{..} \sum f(d) \hphantom{...}$\\ 
$\hphantom{..}1 \quad \hphantom{...}  + 2 \quad + 2 $\\
$\hphantom{..}2 \quad \hphantom{...}  - 1 \quad + 1 $\\
$\mu (3) = - 1 \quad \quad0$

\noindent for $n=4$,\\
$\hphantom{..}d \quad \hphantom{.} f(d) \hphantom{..} \sum f(d) \hphantom{...}$\\ 
$\hphantom{..}1 \quad \hphantom{.}  + 3 \quad + 3 $\\
$\hphantom{..}2 \quad \hphantom{.}  - 2 \quad + 1 $\\
$\hphantom{..}3 \quad \hphantom{.}  - 1 \quad \hphantom{...} 0 $\\
$\mu (4) = \hphantom{.}0 \quad \hphantom{...} 0$

\noindent for $n = 5$,\\
$\hphantom{..}d \quad \hphantom{..}  f(d) \hphantom{..} \sum f(d)$\\
$\hphantom{..}1 \quad \hphantom{..}  + 4 \quad + 4 $\\
$\hphantom{..}2 \quad \hphantom{..}  - 2 \quad + 2 $\\
$\hphantom{..}3 \quad \hphantom{..}  - 1 \quad + 1 $\\
$ \mu (5) = - 1 \quad \hphantom{.. ..} 0$

\noindent for $ n = 6$,\\ 
$\hphantom{..}d \quad \hphantom{..} f(d) \hphantom{..} \sum f(d)$\\
$\hphantom{..}1 \quad \hphantom{..}  + 5 \quad + 5 $\\
$\hphantom{..}2 \quad \hphantom{..}  - 3 \quad + 2 $\\
$\hphantom{..}3 \quad \hphantom{..}  - 2 \quad \quad 0 $\\
$\hphantom{..}5 \quad \hphantom{..}  - 1 \quad - 1 $\\
$ \mu (6) = + 1 \quad \hphantom{...} 0$

\newpage

\noindent On page 217, Spencer-Brown says:
``All it (the cascade) requires is what previous cascades have told it. We also see that there is no need to list the final term for any $n$, since the answer must be the penultimate term in the running total with the sign reversed". Thus, the penultimate term in $\sum f(d)$ with reversed sign
$$= -\bigg(\sum_{d=1}^{n-1} f(d)\bigg) = f(n) \hphantom{... .} (= \text{the final } f(d) \text{ term} = \mu(n)\big).$$
\noindent For Spencer-Brown's cascade in Appendix 9 \cite{SBG}, note that the final term of $\sum f(d) = 0$. In his letter to Moshe Klein, however, Spencer-Brown provides a variation of the cascade where $\sum_{d=1}^{n} f(d) = \mu(n)$. Consider,\\
\noindent $d$ \quad \quad $n = 2$,\\ 
\hphantom{1} \quad $ - 1 \quad - 1 $\\
$ \therefore \ \mu (2) = - 1$

\noindent $d$ \quad \quad $ n = 3$\\
\hphantom{1} \quad $ - 2 \quad - 2 $\\
$2$ \quad $ + 1 \quad - 1 $\\
$ \therefore \ \mu (3) = - 1$

\noindent $d$ \quad \quad $ n = 4$\\ 
\hphantom{1} \quad $ - 3 \quad - 3 $\\
$2$ \quad $ + 2 \quad - 1 $\\
$3$ \quad $ + 1 \quad \quad 0 $\\
$ \therefore \ \mu (4) = 0$

\noindent $d$ \quad \quad $ n = 5$\\ 
\hphantom{1} \quad $ - 4 \quad - 4 $\\
$2$ \quad $ + 2 \quad - 2 $\\
$3$ \quad $ + 1 \quad - 1 $\\
$ \therefore \ \mu (5) = - 1$

\noindent $d$ \quad \quad $ n = 6$\\ 
\hphantom{1} \quad $ - 5 \quad - 5 $\\
$2$ \quad $ + 3 \quad - 2 $\\
$3$ \quad $ + 2 \quad \quad 0 $\\
$5$ \quad $ + 1 \quad + 1 $\\
$ \therefore \ \mu (6) = + 1$

\newpage

\noindent We begin to see how Spencer-Brown rectifies Legendre's method of counting large primes in $n$ and then corrupts it to give the answer zero for all $(n)$. This is illustrated by the cascade in Appendix 9 \cite{SBG}, in which the final term of the ``running total" or $\sum f(d)$ is $0$. 

\noindent For $n=20$, which is a squareful number, we have $f(20)= \mu(20) = 0$ and $\sum_{d=1}^{20} f(d) = 0$. Spencer-Brown says, `` double $n$ to $2n = 40$. The aggregate of pluses and minuses in the upper section of $2n = 40$ is exactly what it was in the whole of $n = 20$. But it contains two more minus signs than did the upper section of $n=20$, so its sum is likely to be reduced towards or beyond zero." Spencer-Brown continues, ``every time we doubled the argument we would have to add an average of two more minus signs to the upper-section terms of the new doubled argument, so there must come a time when the sum of the upper-section terms of the new doubled argument is reduced to or beyond zero. Suppose it is reduced to zero. Then the sum of the lower-section terms for this argument will also be zero, and there will be no tendency in either direction when it is doubled again. But suppose the sum in the upper section is reduced beyond zero to a negative value. Now the lower-section sum for this argument must be positive, and the whole process must play itself out again, this time in the opposite direction" (page 205-206, Appendix 9 \cite{SBG}). 

\noindent Recall that Spencer-Brown splits $n \leq d$ into two sections at $\frac{n}{2}$, calling $d \leq \frac{n}{2}$ the upper section and $\frac{n}{2} < d \leq n$ the lower section. For the lower section, $\mu(d)\floor*{\frac{n}{d}} = \mu(d)$. So if $n=20$, the lower section of $\sum f(d)$ is $\sum_{d=11}^{20} \mu(d) = -2$, and this sum balances the upper section, which is $\sum_{d=1}^{10} \mu(d)\floor*{\frac{n}{d}} = +2$. These two sums equal each other, but are reversed in sign since they sum to $0$. For $2n = 40$, the lower section of $\sum f(d)$ is $\sum_{d=21}^{40} \mu(d) = +3$ and the upper section of $\sum f(d)$ is $\sum_{d=1}^{20} \mu(d)\floor*{\frac{40}{d}} = -3$. For $4n = 80$, the lower section of $\sum f(d)$ is $\sum_{d=41}^{80} \mu(d) = -4$ and the upper section of $\sum f(d)$ is $\sum_{d=1}^{20} \mu(d)\floor*{\frac{80}{d}} = +4$.

\noindent One may especially note that the final $f(d)$ term, $f(n)=\mu(n)\floor*{\frac{n}{n}}=\mu(n)= +1, -1$ or $0$, and that the final $\sum f(d)$ term is always $0$. The cascade in Appendix 9 \cite{SBG} also illustrates that if $\mu(n)$ is the $n^{th}$ term of $f(d)$, and $\mu(n) = \pm 1$, then the $(n-1)^{th}$ term of $\sum f(d) = \mp 1$, for square-free $n$.

%\noindent What makes $n = 20$ important is that it is a squareful number and $\mu(n) = 0$. Thus, for any squareful number $n$, such as $n=20$, we have $f(20)= \mu(20) = 0$ and $\sum_{d=1}^{20} f(d) = 0$.

%Using the fact that for any squareful $n$, 
%$$\sum_{d=1}^{n} f(d) = 0,$$

%\noindent By using the generalizable example of $n = 20$, and doubling $n$ to $2n = 40$, and then doubling $2n$ to $4n = 80$, and so on, the lower sections eventually include all values of $\mu(d)$ for square-free $(d)>1$" (page 214, Appendix 9 \cite{SBG}), and the ``average algebraic sum, i.e., the sum divided by the number of LSB terms displayed, in each section varies around and is asymptotic to zero" (page 214, Appendix 9, \cite{SBG}), since

%Thereupon, Spencer-Brown says that ``any difference between the upper and the lower sections can be seen to be self-annihilating, because of the negative feedback between the two sections" (page 206, \cite{SBG}), which is to say that the lower section and the upper section must sum to zero.

\newpage
\thispagestyle{empty}
\section*{Appendix 0 -- A M\"{o}bius and Liouville Excursion}

\noindent Let $\{ p_1, p_2, p_3, \cdots , p_n\}$ be a given set of prime numbers. Note the identity that 
$$(1+p_1)(1+p_2) \ldots (1+p_n)$$ is equal to  the sum of all possible products of $p_i$ with exactly one or no appearance of each $p_i$. 
This implies directly, by the binomial theorem, that the sum of the M\"{o}bius ($\mu$) values of each of these products  is equal to $(1+(-1))^n = 0$ (here $n$ is the number of primes being considered).
From this we can conclude well-known identities such as $\sum_{d | n} \mu(d) = 0$ where $d$ runs over all the divisors of $n$ including $n$ itself. This appendix provides points of view about this identity in the context of the 
Legendre/Spencer-Brown formula and the Cascade of the previous appendix.

\noindent Recall that the Legendre/Spencer-Brown (LSB) formula\footnote{The Legendre formula, rectified by Spencer-Brown, is presented as follows:
$$ \pi (n) = \pi ( \sqrt{n}) + \sum_{d = 1}^{\sqrt{n}} \mu (d) \left[ \frac{n}{d} \right].$$
Specifically for $d = 1$, Spencer-Brown rectifies the dividend to $\left[\frac{n-1}{d}\right]$. We take the bracketed variable $\left[ \frac{n}{d} \right]$ to mean the floored quotient $\floor*{\frac{n}{d}}$. Further, $d$ is any number $\leq n$ that can be expressed as a product of primes up to $\sqrt{n}$.} requires all primes up to $\sqrt{n}$ to find all primes up to $n$. For the sake of computing $\pi (n)$, build all square-free divisors $d \leq n,$ by starting with the first square-free divisor $1$.

\noindent Call:

\noindent $\{p_i\} = \{p_i | ``p_i$ \text{is prime" and} ``$1 \leq i \leq \pi(\sqrt{n})"\}$\\
$\{p_i.p_j\} = \{(p_i.p_j)|``p_i,p_j$ \text{are prime" and } ``$1 \leq i < j \leq \pi(\sqrt{n})"\}$\\
$\{p_i.p_j.p_k\} = \{(p_i.p_j.p_k)|``p_i, p_j, p_k$ \text{are prime" and } $``1 \leq i < j < k \leq \pi(\sqrt{n})"\}$\\
...

\noindent Note how the small primes up to $\sqrt{n}$ combine uniquely to form a large (but finite) collection of square-free divisors $d = \{1, \{p_i\}, \{p_i.p_j\}, \{p_i.p_j.p_k\}, ... \}$. The last item in this list is the product of all the primes up to $\sqrt{n}$. By definition of the M\"{o}bius function, we know that $\mu (1) = +1$, $\mu(p_i) = -1$, $\mu (p_i.p_j) = +1$, $\mu (p_i.p_j.p_k) = -1$, ..., for all primes $p_i, p_j$, and $p_k$. We wish to keep track of the $+1$'s and $-1$'s contributed by these square-free divisors.

\newpage

\noindent \noindent Let $\#S =$ the number of elements of a set $S$. Let $\#\{p_i\} = m$. For $m = 5$,

\noindent $\#\{p_i\} = \#\{2, 3, 5, 7, 11\} = 5 = m$, in general.

\noindent $\#\{p_i.p_j\} = \#2 \{3, 5, 7, 11\} = \#\{2.3, 2.5, 2.7, 2.11\} = 4 = (m-1)$ elements, in general.\\
\hphantom{+} \quad \quad \quad $+\#3 \{5, 7, 11\}, = \#\{3.5, 3.7, 3.11\} = 3$ elements\\
\hphantom{+} \quad \quad \quad $+\#5 \{7, 11\}, = \#\{5.7, 5.11\} = 2$ elements\\
\hphantom{+}\quad \quad \quad $+\#7 \{11\} = \#\{7.11\} = 1$ element.

\noindent Note that $\{p_i.p_j\}$ has $1+2+3+4$ elements, and in general, $\#\{p_i.p_j\} = 1+2+3+...+(m-1)$.

\noindent $\#\{p_i.p_j.p_k\} = \#2\begin{cases} 3\{5, 7, 11\}&=\#\{2.3.5, 2.3.7, 2.3.11\}=3=(m-2) \text{ elements}\\ 5\{7, 11\}&=\#\{2.5.7, 2.5.11\}=2 \text{ elements}\\ 7 \{11\}&=\#\{2.7.11\}=1 \text{ element} \end{cases}$\\
\hphantom{+}\quad \quad \quad \quad \quad \quad \quad $+\#3\begin{cases} 5\{7, 11\}&=\#\{3.5.7, 3.5.11\}=2 \text{ elements},\\ 7\{11\}&=\#\{3.7.11\}=1 \text{ element} \end{cases}$\\
\hphantom{+}\quad \quad \quad \quad \quad \quad \quad $+\#5 \begin{cases} 7\{11\}&=\#\{5.7.11\}=1 \text{ element}.
\end{cases}$

\noindent Note that $\{p_i.p_j.p_k\}$ has $1 + (1+2) + (1+2+3)$ elements, and in general,\\ $\#\{p_i.p_j.p_k\} = 1 + (1+2) + (1+2+3) + ... + (1+2+3+...+(m-2))$.

\noindent Recall that $\#\{p_i.p_j\} = 1+2+3+...+(m-2)+(m-1)$. Notice how this sum, can be broken up into two parts, $1+2+3+...+(m-2)$, and $(m-1)$. Take the first part. It is equal to the largest collection in $\{p_i.p_j.p_k\}$ with $(1+2+3+...+(m-2))$ elements.

\noindent Similarly, recall that $\#\{p_i\} = m$. We can break this sum into two parts $(m-1)$ and $1$. Take the first part. It equals the second part $(m-1)$ of $\{p_i.p_j\}$. The second part, $1$ of $\{p_i\}$, however, is equal to the first square-free divisor $1$.

\noindent As we collect all square-free products of $p_i, p_j, p_k, ...$ into $\{1, \{p_i\}, \{p_i.p_j\}, \{p_i.p_j.p_k\}, ... \}$, we see the $1$ giving $\mu(1) = +1$, the $\{p_i\}$ giving $\mu (p_i) = -1$, the $\{p_i.p_j\}$ giving $\mu(p_i.p_j) = +1$, and so on, and we observe how one collection of $+1$'s cancels out a preceding or successive collection of $-1$'s in our set of divisors.
\newpage

\noindent $\{1\}$ has $1$ element.\\
\noindent $\{p_i\}$ has $m$ elements = $(m-1)$ $+$ $1$

\noindent $\{p_i.p_j\}$ has \begin{tabular}{c}
%\hline
1+2+3+\ldots+$(m-2)$\\
\end{tabular} +$(m-1)$ elements.

\noindent $\{p_i.p_j.p_k\}$ has
\begin{tabular}{c}
%\hline
1\\
+1+2\\
+1+2+3\\
+...\\
+1+2+3+\ldots+$(m-3)$\\
\end{tabular}
+ \begin{tabular}{c}
%\hline
1+2+3+\ldots+$(m-2)$\\
\end{tabular} elements.

\noindent $\{p_i.p_j.p_k.p_l\}$ has 1
+ \begin{tabular}{c}
    %\hline
    1  \\
    +1+2 \\
\end{tabular}
+ \begin{tabular}{c}
    %\hline
    1\\
    +1+2\\
    +1+2+3\\
\end{tabular}
+ \quad ...
+ \begin{tabular}{c}
    %\hline
    1\\
    +1+2\\
    +1+2+3\\
    +\ldots\\
    +1+2+3+\ldots+$(m-4)$
\end{tabular}\\
\vspace{8 pt}
\begin{center}
+ \begin{tabular}{c}
    %\hline
    1\\
    +1+2\\
    +1+2+3\\
    +\ldots\\
    +1+2+3+\ldots+$(m-3)$\\
\end{tabular} elements.
\end{center}

%\noindent $\{p_i.p_j.p_k.p_l\}$ has 
%\begin{cases} 1\\
%+\begin{cases}1\\+1+2\end{cases}\\
%+\begin{cases}1\\+1+2\\+1+2+3\end{cases}\\
%+ \quad ...\\
%+\begin{cases}1\\+1+2\\+1+2+3\\+\quad...\\+1+2+3+...+$(m-4)$\end{cases}\\
%+\begin{cases}1\\+1+2\\+1+2+3\\+\quad...\\+1+2+3+...+$(m-3)$\end{cases}\\ 
%\end{cases} elements.

\newpage

%\includepdf[scale=0.75]{equiprobable}

\noindent Notice how the total number of divisors in each collection correspond to the counts of the vertical columns of Pascal's triangle:

%>{$}l<{$}|*{7}{c}
\begin{center}
\begin{tabular}{c c c c c c c}
1&&&&&&\\
1&1&&&&&\\
1&2&1&&&&\\
1&3&3&1&&&\\
1&4&6&4&1&&\\
1&5&10&10&5&1&\\
1&6&15&20&15&6&1\\
\end{tabular}
\end{center}

\noindent We see how column $1$ corresponds to the total number of divisors in $\{p_i\}$,\\
column $2$ corresponds to the total number of divisors in $\{p_i.p_j\}$,\\
column $3$ corresponds to the total number of divisors in $\{p_i.p_j.p_k\}$,\\
column $4$ corresponds to the total number of divisors in $\{p_i.p_j.p_k.p_l\}$,\\
...

\noindent \textbf{Column $\mathbf{1}$} $= 1+1+1+1+1+1+1+... = \#\{p_i\}$ with $\mu(d) = -1$.\\
\textbf{Column $\mathbf{2}$} $= 1+2+3+4+5+6+ ... = \#\{p_i.p_j\}$ with $\mu(d) = +1$. If we look at the $6$ in this list, it is the sum of all the entries of column $1$ till $5$, i.e., $1+1+1+1+1+1=6$.\\
\textbf{Column $\mathbf{3}$} is the sum of \textit{triangular} numbers $= 1+3+6+10+15+ ... = \#\{p_i.p_j.p_k\}$ with $\mu(d) = -1$. If we look at the $15$ in this list, it is the sum of all entries of column $2$ till $10$, i.e., $1+2+3+4+5=15$.\\ 
\textbf{Column $\mathbf{4}$} is the sum of \textit{tetrahedral} numbers $= 1+4+10+20+35+ ... =\#\{p_i.p_j.p_k.p_l\}$ with $\mu(d) = +1$. If we look at the $20$ in this list, it is the sum of all entries of column $3$ till $10$, i.e., $1+3+6+10=20$.\\
\textbf{Column $\mathbf{5}$} is the sum of \textit{pentatope} numbers $= 1+5+15+35+70+ ... = \#\{p_i.p_j.p_k.p_l.p_m\}$ with $\mu(d) = -1$. If we look at the $15$ in this list, it is the sum of all entries of column $4$ till $5$, i.e., $1+4+10=15$. 

\newpage

\noindent Note: Pascal referred to this self reference as a ``hockey stick pattern". When the numbers are arranged in Pascal's form, the rows and columns are interchangeable as follows: \begin{center}
\begin{tabular}{c c c c c c c}
1&1&1&1&1&1&1\\
1&2&3&4&5&6&...\\
1&3&6&10&15&...&\\
1&4&10&20&...&&\\
1&5&15&...&&&\\
1&6&...&&&&\\
1&...&&&&&\\%\hline
\end{tabular}
\end{center}
The ``hockey stick pattern" that Pascal is referring to is identical to the dovetailing identities we saw in the previous page. If we look at the $6$ in Column $2$, it is the sum of all numbers in the previous Column $1$ up to the row of $6$ itself. This gives us the identity $1+1+1+1+1+1=6$, which is the same as what we had on the previous page. Similarly, if we look at the $15$ in Column $3$, it is the sum of all numbers in the previous Column $2$ up to the row of $15$ itself. This gives us the identity $1+2+3+4+5=15$, which is the same as what we had on the previous page. Similarly, if we look at the $20$ in Column $4$, it is the sum of all numbers in the previous Column $3$ up to the row of $20$ itself. This gives us the identity $1+3+6+10=20$. The $15$ in Column $5$ is the sum of all numbers in the previous Column $4$ up to the row of $15$ itself, giving us the identity $1+4+10=15$. The $6$ in Column $6$ is the sum of all numbers in the previous Column $5$ up to the row of $6$ itself, giving us the identity $1+5=6$. Finally, the $1$ in the last Column $7$ equals the adjacent $1$ in the previous Column $6$:

\noindent $1+1+1+1+1+1=6$\\
$1+2+3+4+5=15$\\
$1+3+6+10=20$\\
$1+4+10=15$\\
$1+5=6$\\
$1=1$

\newpage

\noindent Pascal said\footnote{\href{https://www.cut-the-knot.org/arithmetic/combinatorics/PascalTriangleProperties.shtml}{https://www.cut-the-knot.org/arithmetic/combinatorics/PascalTriangleProperties.shtml}},``In every arithmetical triangle each cell is equal to the sum of all the cells of the preceding row from its column to the first, inclusive". Since the rows and columns are interchangeable, we also have ``In every arithmetical triangle each cell is equal to the sum of all the cells of the preceding column from its row to the first, inclusive".

\noindent In collecting all square-free divisors from a set of prime numbers $\{p_i\}$ of size $m$, we notice that the count of all possible combinations of square-free divisors generated by the set $\{p_i\}$, when collected combinatorially in collections of $\{p_i\}, \{p_i.p_j\}$, $\{p_i.p_j.p_k\}, \ldots $ , correspond exactly to the counts in the columns of Pascal's triangle. If the columns of Pascal's triangle are assigned $-$ or $+$ in an alternating way, i.e., if Column $1$ is assigned $-$, Column $2$ is assigned $+$, Column $3$ is $-$, Column $4$ is $+$,..., then not only do these columns correspond to the counts of square-free divisors for each collection $\{p_i\}, \{p_i.p_j\}, \{p_i.p_j.p_k\}, ...$, but also each column now corresponds to the cumulative sum of M\"{o}bius values of the square-free divisors that correspond to that column, contributing a positive M\"{o}bius sum for columns $2, 4, 6, 8, ...$ or negative M\"{o}bius sum for columns $1, 3, 5, 7, ...$ .

\noindent Column 1, which now represents the singular prime count up to $m$, gives $\mu(d) = -1$ throughout the column. So, if we sum all members of Column 1, we are summing the M\"{o}bius values of each of the divisors that Column 1 corresponds to, which is simply $\#\{p_i\}$.  Column 2, which represents all possible $(p_i.p_j)$ combinations, gives $\mu(d) = +1$. So, if we sum all members of Column $2$, we are summing the M\"{o}bius values of each of the divisors that Column $2$ corresponds to, which is simply $\#\{p_i.p_j\}$. Column $3$, which represents all possible $(p_i.p_j.p_k)$ combinations, gives $\mu(d) = -1$. So, if we sum all members of Column $3$, we are summing the M\"{o}bius values of each of the divisors that Column $3$ corresponds to, which is simply $\#\{p_i.p_j.p_k\}$.

\noindent Pascal's ``hockey stick pattern" shows us that the alternating sum of all columns of a finite Pascal's triangle equals -1. One may verify this for a Pascal's triangle of any size. For a Pascal's triangle of size $m = 7$, the alternating sum of all $7$ columns\\
$= -(1+1+1+1+1+1+1)+(1+2+3+4+5+6)-(1+3+6+10+15)+(1+4+10+20)-(1+5+15)+(1+6)-1 = \mathbf{-1}.$
%$= -(6+1)+(15+6)-(20+15)+(15+20)-(6+15)+(1+6)-1$,\\
%$= -6-1+15+6-20-15+15+20-6-15+1+6-1,$\\
%$= -1$.\\

\noindent Column $1$ has one remaining element $p_*$ that survives as all the other columns cancel each other out in the alternate sign summation of all columns of Pascal's triangle. We know $\mu(p_*)=-1$. We see that $p_*$ could be none other than the first prime $2$, as shown in the first Pascal arrangement, or the $m$-th prime as seen in Pascal's original arrangement\footnote{\href{http://oeis.org/wiki/Pascal_triangle}{http://oeis.org/wiki/Pascal\_triangle}}. This taken together with $\mu(1)=+1$ will always give us a total M\"{o}bius sum of zero for all square-free divisors generated by a finite set of primes, including the set of primes and the first square-free divisor $1$.
%http://oeis.org/wiki/Pascal_triangle#Pascal.27s_triangle_rows_alternating_sign_sums
\begin{thm}{\hphantom{.}}
A. Let $n$ be a natural number. Let $p_1, p_2, p_3, ...,$ enumerate the primes up to $\sqrt{n}$. Let $d$ be a square-free divisor in $\{1, \{p_i\}, \{p_i.p_j\}, \{p_i.p_j.p_k\}, ... \}$, where\\
\indent $\{p_i\}:= \{p_i | ``p_i$ is prime" and ``$1 \leq i \leq \pi(\sqrt{n})"\}$,\\
\indent $\{p_i.p_j\}:= \{(p_i.p_j)|``p_i,p_j$ are prime" and ``$1 \leq i < j \leq \pi(\sqrt{n})"\}$\\
\indent $\{p_i.p_j.p_k\}:= \{(p_i.p_j.p_k)|``p_i, p_j, p_k$ are prime" and $``1 \leq i < j < k \leq \pi(\sqrt{n})"\}$, ...\\ Then, $$\sum_{d}^{} \mu (d) = 0.$$\\
B. Let $n$ be a natural number. Let $p_1, p_2, p_3, ..., p_n$ be any collection of $n$ primes.\\ 
Let $d$ be a square-free divisor in $\{1, \{p_i\}, \{p_i.p_j\}, \{p_i.p_j.p_k\}, ... \}$, where\\
\indent $\{p_i\}:= \{p_i | ``p_i$ is prime" and ``$1 \leq i \leq n"\}$,\\
\indent $\{p_i.p_j\}:= \{(p_i.p_j)|``p_i,p_j$ are prime" and ``$1 \leq i < j \leq n"\}$\\
\indent $\{p_i.p_j.p_k\}:= \{(p_i.p_j.p_k)|``p_i, p_j, p_k$ are prime" and $``1 \leq i < j < k \leq n"\}$, ...\\ Then, $$\sum_{d}^{}\mu (d) = 0.$$
\end{thm}
\noindent We see how any collection of $n$ primes can form divisors $d$ in $\{1, \{p_i\}, \{p_i.p_j\}$, $\{p_i.p_j.p_k\}, ... \}$ and that $\#\{p_i\}, \#\{p_i.p_j\}, \#\{p_i.p_j.p_k\}, ...$ correspond exactly to the columns (or rows) of Pascal's triangle. The ``alternating sums of subsequent columns (or rows) in Pascal's triangle" is $\{1,0,0,0,,... \}$, which gives a total M\"{o}bius sum of -1. Since $\mu (1) = +1$, $$\sum_{d}^{}\mu(d)=\mu(1)-1=0.$$

\newpage

\noindent The Liouville function $\lambda$ maps to values $+1$ or $-1$. If the sum of the powers of the prime components of a number $n$ is even, $\lambda (n) = +1$, and if the sum of the powers of the prime components of $n$ is odd, $\lambda (n) = -1$.

\noindent Consider, the set of primes $\{2, 3\}$. Write the ``Pascal form" or ``triangular form" of the collection of divisors of the maximal product of elements in that set, as: $\begin{pmatrix}
  2 & 3\\ 
  \phantom & 2.3
\end{pmatrix}$. This is the collection of all square-free products generated by $\{2, 3\}$.

Now consider, $\begin{pmatrix}
  \mathbf{2}(2) & \mathbf{3}(3)\\ 
  \phantom & \mathbf{2.3} \begin{pmatrix}
  2 & 3\\ 
  \phantom & 2.3
\end{pmatrix}
\end{pmatrix}$ 
$= \begin{pmatrix}
  2^2 & 3^2\\ 
  \phantom & \begin{pmatrix}
  2^2.3 & 2.3^2\\ 
  \phantom & 2^2.3^2
\end{pmatrix}
\end{pmatrix}$.

\noindent We see how a simple product of each entry of this Pascal form (with itself or its corresponding Pascal form), helps us generate all the possible combinations of $\{2,3\}$ of degree $2$. Call the degree\footnote{This type of definition has a similar sense to the ``degree" of a polynomial.} of a number the degree of the highest-degree-prime in the number's prime factorization. For instance, the degree of $45 = 3^2.5$ is $2$.

\noindent Now consider Degree $3$ collections:
$$\begin{pmatrix}
  2(2(2)) & 3(3(3))\\ 
  \phantom & 2.3 \begin{pmatrix}
  2(2) & 3(3)\\ 
  \phantom & 2.3 \begin{pmatrix}
  2 & 3\\
  \phantom & 2.3 \end{pmatrix}
\end{pmatrix}
\end{pmatrix}
= \begin{pmatrix}
  2^3 & 3^3\\ 
  \phantom  &\begin{pmatrix}
  2^3.3 & 2.3^3\\ 
  \phantom &\begin{pmatrix}
  2^3.3^2 & 2^2.3^3\\
  \phantom & 2^3.3^3 \end{pmatrix}
\end{pmatrix}
\end{pmatrix}$$
\noindent Degree $4$ collections:
$$\begin{pmatrix}
  2(2(2(2))) & 3(3(3(3)))\\ 
  \phantom & 2.3 \begin{pmatrix}
  2(2(2)) & 3(3(3))\\ 
  \phantom & 2.3 \begin{pmatrix}
  2(2) & 3(3)\\
  \phantom & 2.3 \begin{pmatrix}
  2 & 3\\
  \phantom & 2.3 \end{pmatrix} \end{pmatrix}
\end{pmatrix}
\end{pmatrix}$$
$$ = \begin{pmatrix}
  2^4 & 3^4\\ 
  \phantom & 2.3 \begin{pmatrix}
  2^3 & 3^3\\ 
  \phantom & 2.3 \begin{pmatrix}
  2^2 & 3^2\\
  \phantom & 2.3 \begin{pmatrix}
  2 & 3\\
  \phantom & 2.3 \end{pmatrix} \end{pmatrix}
\end{pmatrix}
\end{pmatrix}
$$
$$ = \begin{pmatrix}
  2^4 & 3^4\\ 
  \phantom & \begin{pmatrix}
  2^4.3 & 2.3^4\\ 
  \phantom & \begin{pmatrix}
  2^4.3^2 & 2^2.3^4\\
  \phantom & \begin{pmatrix}
  2^4.3^3 & 2^3.3^4\\
  \phantom & 2^4.3^4 \end{pmatrix} \end{pmatrix}
\end{pmatrix}
\end{pmatrix}
$$
\noindent Degree $5$ collections:
$$\begin{pmatrix}
  2(2(2(2(2)))) & 3(3(3(3(3))))\\ 
  \phantom & 2.3 \begin{pmatrix}
  2(2(2(2))) & 3(3(3(3)))\\ 
  \phantom & 2.3 \begin{pmatrix}
  2(2(2)) & 3(3(3))\\
  \phantom & 2.3 \begin{pmatrix}
  2(2) & 3(2)\\
  \phantom & 2.3 \begin{pmatrix}
  2 & 3\\
  \phantom & 2.3 \end{pmatrix} \end{pmatrix} \end{pmatrix}
\end{pmatrix}
\end{pmatrix}
$$
$$ = \begin{pmatrix}
  2^5 & 3^5\\ 
  \phantom & 2.3 \begin{pmatrix}
  2^4 & 3^4\\ 
  \phantom & 2.3 \begin{pmatrix}
  2^3 & 3^3\\
  \phantom & 2.3 \begin{pmatrix}
  2^2 & 3^2\\
  \phantom & 2.3 \begin{pmatrix}
  2 & 3\\
  \phantom & 2.3 \end{pmatrix} \end{pmatrix} \end{pmatrix}
\end{pmatrix}
\end{pmatrix}
$$
$$ = \begin{pmatrix}
  2^5 & 3^5\\ 
  \phantom & \begin{pmatrix}
  2^5.3 & 2.3^5\\ 
  \phantom & \begin{pmatrix}
  2^5.3^2 & 2^2.3^5\\
  \phantom & \begin{pmatrix}
  2^5.3^3 & 2^3.3^5\\
  \phantom & \begin{pmatrix}
  2^5.3^4 & 2^4.3^5\\
  \phantom & 2^5.3^5 \end{pmatrix} \end{pmatrix} \end{pmatrix}
\end{pmatrix}
\end{pmatrix}
$$

\noindent If we were to replace every divisor in these collections with $o$ for ``odd", when $\lambda (d) = -1$ and $e$ for ``even", when $\lambda (d) = +1$, we have:
\newpage

\noindent Degree $0$ : \quad $e$\\
Degree $1$: $\begin{pmatrix}
  o & o\\ 
  \phantom & e
\end{pmatrix}$\\
Degree $2$: $\begin{pmatrix}
  e & e\\ 
  \phantom & \begin{pmatrix}
  o & o\\ 
  \phantom & e
\end{pmatrix}
\end{pmatrix}$\\
Degree $3$: $\begin{pmatrix}
  o & o\\ 
  \phantom  &\begin{pmatrix}
  e & e\\ 
  \phantom &\begin{pmatrix}
  o & o\\
  \phantom & e \end{pmatrix}
\end{pmatrix}
\end{pmatrix}
$\\
Degree $4$: $\begin{pmatrix}
  e & e\\ 
  \phantom & \begin{pmatrix}
  o & o\\ 
  \phantom & \begin{pmatrix}
  e & e\\
  \phantom & \begin{pmatrix}
  o & o\\
  \phantom & e \end{pmatrix} \end{pmatrix}
\end{pmatrix}
\end{pmatrix}
$\\
and Degree $5$: $\begin{pmatrix}
  o & o\\ 
  \phantom & \begin{pmatrix}
  e & e\\ 
  \phantom & \begin{pmatrix}
  o & o\\
  \phantom & \begin{pmatrix}
  e & e\\
  \phantom & \begin{pmatrix}
  o & o\\
  \phantom & e \end{pmatrix} \end{pmatrix} \end{pmatrix}
\end{pmatrix}
\end{pmatrix}$

\noindent We wish to sum the $\lambda$ values over these collections:

\noindent Degree $0$ is simply $e$ or $+1$, and Degree $1$, by Theorem $1$, sums to $o$ or $-1$. Degree $2$ sums to $e$, Degree $3$ sums to $o$, Degree $4$ sums to $e$, Degree $5$ sums to $o$, ...

\noindent Consider, now, the set of primes $\{2, 3, 5\}$. The degree $1$ divisors are $$\begin{pmatrix} 2 & 3 & 5\\
\phantom & 2.3 & 2.5\\
\phantom & \phantom & 3.5\\ 
\phantom & \phantom & 2.3.5 \end{pmatrix}$$,

\newpage

\noindent the degree $2$ divisors are: $\begin{pmatrix} 2^2 & 3^2 & 5^2\\
\phantom & 2.3 \begin{pmatrix} 2 & 3\\
\phantom & 2.3 \end{pmatrix} & 2.5 \begin{pmatrix} 2 & 5\\
\phantom & 2.5 \end{pmatrix}\\
\phantom & \phantom & 3.5 \begin{pmatrix} 3 & 5\\
\phantom & 3.5 \end{pmatrix}\\ 
\phantom & \phantom & 2.3.5 \begin{pmatrix} 2 & 3 & 5\\
\phantom & 2.3 & 2.5\\
\phantom & \phantom & 3.5\\
\phantom & \phantom & 2.3.5 \end{pmatrix} \end{pmatrix}$,

\noindent and the degree $3$ divisors are:
$$\begin{pmatrix} 2^3 & 3^3 & 5^3\\
\phantom & 2.3 \begin{pmatrix} 2^2 & 3^2\\
\phantom & 2.3 \begin{pmatrix} 2 & 3\\
\phantom & 2.3 \end{pmatrix} \end{pmatrix} & 2.5 \begin{pmatrix} 2^2 & 5^2\\
\phantom & 2.5 \begin{pmatrix} 2 & 5\\
\phantom & 2.5 \end{pmatrix} \end{pmatrix}\\
\phantom & \phantom & 3.5 \begin{pmatrix} 3^2 & 5^2\\
\phantom & 3.5 \begin{pmatrix} 3 & 5\\
\phantom & 3.5 \end{pmatrix} \end{pmatrix}\\ 
\phantom & \phantom & 2.3.5 \begin{pmatrix} 2^2 & 3^2 & 5^2\\
\phantom & 2.3 \begin{pmatrix} 2 & 3\\
\phantom & 2.3 \end{pmatrix} & 2.5 \begin{pmatrix} 2 & 5\\
\phantom & 2.5 \end{pmatrix}\\
\phantom & \phantom & 3.5 \begin{pmatrix} 3 & 5\\
\phantom & 3.5 \end{pmatrix}\\
\phantom & \phantom & 2.3.5 \begin{pmatrix} 2 & 3 & 5\\
\phantom & 2.3 & 2.5\\
\phantom & \phantom & 3.5\\ 
\phantom & \phantom & 2.3.5 \end{pmatrix} \end{pmatrix} \end{pmatrix}$$

\noindent If we were to replace every entry in these collections with $o$ for ``odd", when $\lambda (d) = -1$ and $e$ for ``even", when $\lambda (d) = +1$, we have:

\newpage

\noindent Degree $0$ : \quad $e$\\
Degree $1$: $\begin{pmatrix} o & o & o\\
\phantom & e & e\\
\phantom & \phantom & e\\ 
\phantom & \phantom & o \end{pmatrix}$\\
Degree $2$: $\begin{pmatrix} e & e & e\\
\phantom & \begin{pmatrix} o & o\\
\phantom & e \end{pmatrix} & \begin{pmatrix} o & o\\
\phantom & e \end{pmatrix}\\
\phantom & \phantom & \begin{pmatrix} o & o\\
\phantom & e \end{pmatrix}\\ 
\phantom & \phantom & \begin{pmatrix} e & e & e\\
\phantom & o & o\\
\phantom & \phantom & o\\
\phantom & \phantom & e \end{pmatrix} \end{pmatrix}$\\
Degree $3$: $\begin{pmatrix} o & o & o\\
\phantom & \begin{pmatrix} e & e\\
\phantom & \begin{pmatrix} o & o\\
\phantom & e \end{pmatrix} \end{pmatrix} & \begin{pmatrix} e & e\\
\phantom & \begin{pmatrix} o & o\\
\phantom & e \end{pmatrix} \end{pmatrix}\\
\phantom & \phantom & \begin{pmatrix} e & e\\
\phantom & \begin{pmatrix} o & o\\
\phantom & e \end{pmatrix} \end{pmatrix}\\ 
\phantom & \phantom & \begin{pmatrix} o & o & o\\
\phantom & \begin{pmatrix} e & e\\
\phantom & o \end{pmatrix} & \begin{pmatrix} e & e\\
\phantom & o \end{pmatrix}\\
\phantom & \phantom & \begin{pmatrix} e & e\\
\phantom & o \end{pmatrix}\\
\phantom & \phantom & \begin{pmatrix} o & o & o\\
\phantom & e & e\\
\phantom & \phantom & e\\ 
\phantom & \phantom & o \end{pmatrix} \end{pmatrix} \end{pmatrix}$

\noindent We wish to sum the $\lambda$ values over these collections: Degree $0$ is just simply $e$ or $+1$, and Degree $1$, by Theorem $1$, sums to $o$ or $-1$. Degree $2$ sums to $e$, Degree $3$ sums to $o$, ...

\newpage

\noindent The degree\footnote{As a reminder, this sense of ``degree" is similar in the context of polynomials} of a number is defined as the degree of the highest-degree-prime in the number's prime factorization. We observe that the parity of the degree of a number satisfies the following theorem\footnote{The two theorems mentioned in Appendix 0 are original and are attributed to the co-authors of this paper, J. M. Flagg, Louis H. Kauffman, and Divyamaan Sahoo.}:

\begin{thm}{\hphantom{.}}
Let $\{p_i\}$ be the finite set of primes of size $n$, with $i \leq n$. Define the ``degree of a number" $k$ as the degree of the highest-degree-prime in the number's prime factorization. Call $d_0$ a degree $0$ number, $d_1$ a degree $1$ number, and $d_k$ a degree $k$ number, formed only by primes from $\{p_i\}$. Then,
$$\text{A. }  \sum_{d_1}^{}\lambda(d) = -\lambda(d_0).$$
$$\text{B. } \sum_{d_k}^{}\lambda(d)= (-1)^k.$$
\end{thm}

\noindent (A) Begin with the finite set of primes of size $n$. The Degree $0$ collection always consists of a single element, the number $1$, which has $0$ number of prime divisors, and since $0$ is an even number, $\lambda(1) = +1$. The Degree $1$ collection corresponds to all the square-free numbers generated by our initial set of generators, which, in this case, is the set of the first $n$ primes. Since we know that all possible distinct combinations of a set of primes taken distinctly can be given by the Pascal or ``triangular" form generated by the set, and that $\lambda(d) = \mu(d)$ for each $d$ in the degree $0$ and degree $1$ collections, so by Theorem $1$, 
 $$\sum_{d_1}^{}\lambda(d) = -\lambda(d_0).$$
 
\noindent (B) Now, let us consider the collection of divisors of degree $k$. In our concrete demonstrations for prime generators $\{2,3\}$ and $\{2,3,5\}$, we offered a recursive method to list all possible combinations of these prime generators for \textit{any} degree. The method employed is simple. In order to construct the collection of all possible divisors of degree $n$ generated by a set of primes, we must first write out the degree $1$ collection, which is simply the Pascal or ``triangular" form of the collection of square-free numbers generated by that given set of primes. Then, we multiply each entry of this representation with its corresponding Pascal form, i.e., we multiply each entry $p_i$ with its Pascal form $p_i$, we multiply each entry $p_i.p_j$ with $\begin{pmatrix}
  p_i & p_j\\ 
  \phantom & p_i.p_j
\end{pmatrix}$, we multiply $p_i.p_j.p_k$ with $\begin{pmatrix}
  p_i & p_j & p_k\\ 
  \phantom & p_i.p_j & p_i.p_k\\
  \phantom & \phantom & p_j.p_k\\
  \phantom & \phantom & p_i.p_j.p_k
\end{pmatrix}$, and so on.

\noindent This recursive method, applied \textit{once}, gives us all possible degree $2$ numbers whose prime components are restricted by the starting set of primes of size $n$, i.e., the degree $2$ collection. When this recursive method is applied \textit{twice}, we get the degree $3$ collection, and when applied $(k-1)$ \textit{times}, we get the degree $k$ collection.

\noindent Note that the degree $k$ collection will have $k$ containers\footnote{The reader is encouraged to make note of this in the concrete demonstrations provided for example and perform investigations for large $n$ on a large piece of paper.} by the end of this recursive method, i.e., the depth of degree $k$ arrangement is $k$. Further, by design, we see that a degree $k$ collection is essentially a degree $1$ collection of degree $(k-1)$ collections, since each entry of the degree $1$ collection combines with its corresponding Pascal form and begins nesting within itself $(k-1)$ times in order to exhaustively list the entire degree $k$ collection. In essence, a degree $k$ collection is a degree $1$ collection, and we know, by Theorem $2$.A, $\sum_{d_1}^{}\lambda(d) = -1$, regardless of the parity of $n$.

\noindent Starting from the innermost depth, we will be able to replace degree $1$ collections with\\ $\sum_{d_1}^{}\lambda(d) = -1$, since we wish to sum all the $\lambda$ values in the degree $k$ collection. As we move inside-out, we subsequently remove containers and simplify the arrangement, and are finally left with a degree $1$ arrangement with $n$ entries in the first row, namely, the single primes raised to the $k^{th}$ power. From our knowledge of the Pascal triangle, only one entry of the first row will survive the alternate sign summation of the columns of Pascal's triangle, which contributes $\lambda$ of $+1$ if $k$ is even and $-1$ if $k$ is odd. Hence,

$$\sum_{d_k}^{}\lambda(d) = (-1)^k$$

%\include{appendix1}
%\include{appendix2}
%\include{appendix3}

%\include{bib}

% The convention here is to use for each label the author's initials
% and the year of publication, but you may feel free to use whatever
% is easiest for you to remember.

%\include{ack}
\end{document}